%% file: main.tex
\documentclass[preprint,12pt]{elsarticle}
\pdfoutput=1

\usepackage{amssymb}
\usepackage{graphicx}
\usepackage{latexsym}
\usepackage{amsmath}
\usepackage{mathtools}
\usepackage{tabularx}
\usepackage{booktabs}
\usepackage{siunitx}
\usepackage[autostyle]{csquotes}
\usepackage[linesnumbered,ruled]{algorithm2e}
\newcommand*\diff{\mathop{}\!\mathrm{d}}
\DeclareMathOperator\erf{erf}
\usepackage{url}
\usepackage{xcolor}
\usepackage{subcaption}
\usepackage{placeins}
\usepackage{hyperref}

\journal{arXiv}

\begin{document}

\begin{frontmatter}


\title{A level-set based space-time finite element approach to the modelling of solidification and melting processes}

\author[1]{L.\ Boledi}
\author[1]{B.\ Terschanski}
\author[2,3]{S.\ Elgeti}
\author[1,4]{J.\ Kowalski}

\address[1]{Aachen Institute for Advanced Study in Computational Engineering Science (AICES), RWTH Aachen University, 52056 Aachen, Germany}
\address[2]{Institute of Lightweight Design and Structural Biomechanics (ILSB), TU Wien, 1040 Vienna, Austria}
\address[3]{Chair for Computational Analysis of Technical Systems (CATS), RWTH Aachen University, 52056 Aachen, Germany}
\address[4]{Abteilung Computational Geoscience, Georg-August-Universit\"at G\"ottingen, 37077 G\"ottingen, Germany}

\begin{abstract}
We present a strategy for the numerical solution of convection-coupled phase-transition problems, with focus on solidification and melting. We solve for the temperature and flow fields over time. The position of the phase-change interface is tracked with a level-set method, which requires knowledge of the heat-flux discontinuity at the interface. In order to compute the heat-flux jump, we build upon the ghost-cell approach and extend it to the space-time finite element method. This technique does not require a local enrichment of the basis functions, such as methods like extended finite elements, and it can be easily implemented in already existing finite element codes. Verification cases for the 1D Stefan problem and the lid-driven cavity melting problem are provided. Furthermore, we show a more elaborate 2D case in view of complex applications.
\end{abstract}


\end{frontmatter}

\input{section1} 

\input{section1_and_a_half}

\input{section2}

\input{section3}

\input{section4}

\input{section5}

\section*{Acknowledgments}
The authors were supported by the Helmholtz Graduate School for Data Science in Life, Earth and Energy (HDS-LEE).
The work was furthermore supported by the Federal Ministry of Economic Affairs and Energy, on the basis of a decision by the German
Bundestag (50 NA 1908). The authors gratefully acknowledge the computing time granted by the JARA Vergabegremium and provided on the JARA Partition part of the supercomputer JURECA at Forschungszentrum Jülich \cite{jureca2018}.

\bibliographystyle{elsarticle-num}
\bibliography{main}

\end{document}

%% file: section1.tex
\section{Introduction}
\label{section1}

Phase-transition processes are important for many engineering and scientific applications.
The driving application for this work is cryosphere physics, e.g.\ for assessing processes at the ice-ocean boundary layer or for the model based development of thermal ice exploration robots. The underlying physical processes are complex and their modelling results in coupled systems of partial differential equations. Thus, efficient and robust numerical methods are needed. An overview of commonly used numerical methods for the description of solidification and melting can be found in \cite{Hu1996}.

In this work, we focus on the convection-coupled phase-change from solid to liquid and vice versa. The material is  assumed to be incompressible within the two phases, but we account for density and thermal conductivity changes across the phase-change interface (PCI). To model the spatio-temporal phase-change process, we need to solve for the flow and temperature fields for both phases over time. The main challenge then lies in the evolution of the PCI. Various numerical schemes are available, of which a brief review can be found in \cite{Maitre2006,elgeti2016}. In general, we can distinguish between two approaches. Interface tracking methods provide an explicit description of the PCI throughout the simulation. The marker and cell method (MAC), for instance, uses a set of marker points that are transported by the fluid \cite{Harlow1965}, but it is computationally inefficient due to the addition of the markers. One application of explicit interface tracking in boiling flows is documented in \cite{Juric1998}. Instead, interface capturing methods implicitly represent the interface with the value of a scalar function \cite{Quarteroni2010}. In this context we mention the volume of fluid approach, where the volume fraction contained in a discrete element is tracked \cite{Welch2000,Katopodes2019.1}. This method has lower storage requirements than the MAC approach and is very common in combination with finite volume discretizations \cite{elgeti2016}, but it presents disadvantages when evaluating the geometry of the interface \cite{Katopodes2019.2}. A widely utilized alternative is the level-set method, where the interface is defined as the zero level set of a continuous pseudo-density function \cite{Osher1988}. Since the level-set function is continuous, this method reduces the difficulties associated with a discontinuous volume fraction and allows to keep the interface sharp \cite{Quarteroni2010,elgeti2016}. In the context of this work we describe the PCI with the level-set method. The resulting level-set function is advected according to the propagation speed of the PCI. Such velocity field depends on local energy conservation across the interface and can be modelled as the Stefan condition \cite{Stefan1891,Chen1997}. This formulation requires us to approximate the heat-flux discontinuity across the interface based on the evolving temperature and velocity fields.

The choice of the discretization scheme plays an important role. Finite difference schemes offer an easy implementation and many authors employ them to solve phase-change problems \cite{Chen1997,Gibou2005}. Alternatively, finite element methods (FEM) provide an increased versatility in terms of domain geometry and are widely used in engineering applications \cite{Gross2006,Valance2009}. However, standard Galerkin FEM fail to capture discontinuous derivatives across the PCI that are needed to evaluate the position of the interface. This issue is addressed by extended finite element methods (XFEM). While originally introduced as an alternative to remeshing in mechanical problems such as crack propagation \cite{Dolbow1999}, the method has also been employed to capture discontinuous gradients in Stefan problems \cite{Chessa2002}. A brief introduction to XFEM in this context can be found in Chapter 4.3 of \cite{Herzog2010}. The idea is to locally enrich the FEM basis with functions that have discontinuous derivatives. The drawback is that one needs to update the nodes that are enriched based on the location of the interface. This changes the number of degrees of freedom over time and subsequently requires repeated reallocation of the finite element system matrices. Thus, we choose a conceptually simpler method of recovering flux discontinuities with standard choices of FEM basis functions. We build upon existing work from Gibou and Fedkiw on the ghost-cell approach to describe an extension of the ghost split to arbitrary FEM meshes \cite{Gibou2002,Gibou2005,Gibou2007}. This method circumvents the need for an adaptive enrichment and it is easier to implement in already existing FEM codes. In particular, we show the ghost-cell method applied to our space-time finite element solver \cite{Behr1992,Behr1992.2}.

This paper is structured as follows: In Section \ref{section2}, we describe the general physical setting and the proposed numerical approach. The governing equations for flow and temperature and their discretization within a space-time FEM framework are introduced in Section \ref{sec:modelling}. In Section \ref{section3}, we describe the level-set method to track the interface and the flux reconstruction algorithm based on the ghost-cell approach. In Section \ref{section4}, we validate our method against the analytical solution of a 1D Stefan problem on structured and unstructured meshes. Then, we show a 2D lid-driven cavity problem with temperature and flow field coupling. The last simulation covers the corner flow around a more elaborate geometry. Finally, in Section \ref{section5} we summarize our results and provide an outlook.

%% file: section1_and_a_half.tex
\section{General approach}
\label{section2}
In this section we introduce the physical setting of the problem and give an overview of the proposed numerical strategy.

\begin{figure}
\centering
	\includegraphics[trim={0.0cm 0.0cm 0.0cm 0.0cm},clip, width=0.35\textwidth]{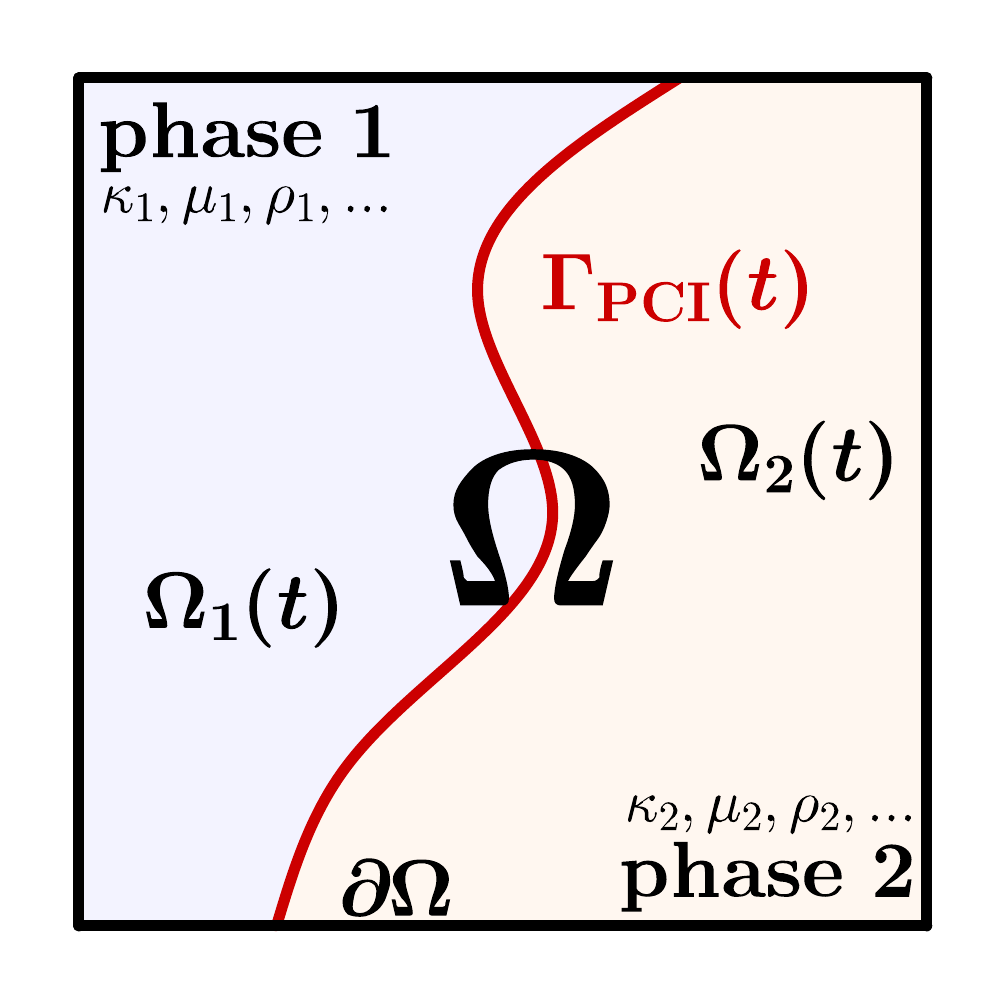}
	\caption{Sketch of the physical situation. The overall domain $\Omega$ consists of subsets $\Omega_1(t)$ and $\Omega_2(t)$ denoting the liquid and solid regions respectively. The red line represents the phase-change interface (PCI). Both the subsets and the PCI change with time, while the domain $\Omega$ and the outer boundary $\partial\Omega$ are constant.}
\label{fig:basicDomain}
\end{figure}

\subsection{Physical setting}
We consider a homogeneous material in a two-phase state: A general domain of interest $\Omega$ consists of a solid region and a liquid region, see Figure \ref{fig:basicDomain}. The objective of our model is to determine the spatio-temporal evolution of the material's temperature and velocity fields, as well as the corresponding evolution of liquid and solid regions. We are mostly interested in complex phase-change processes in water-ice systems, and therefore assume the two phases to be separated from each other by a distinct and well-defined phase-change interface (PCI). Note that this excludes certain alloys or multi-component liquids as well as solutions such as salt water, which tend to develop a transitioning phase-change area rather than a distinct PCI, the so-called mushy layer \cite{worster1997}. The solid and the liquid phase, denoted by the subscript $i\in\{1,2\}$ in Figure \ref{fig:basicDomain}, are each assigned phase-wise constant material parameters for thermal conductivities $\kappa_i$, viscosities $\mu_i$ and densities $\rho_i$. Although this still is an idealization with respect to reality, in which material properties might furthermore obey a continuous temperature sensitivity, e.g. temperature-sensitive density of water and ice \cite{Ulamec2007}, it is a significant step towards improving the predictability of complex phase-change simulations.

\subsection{Overview of the numerical approach}
\begin{figure}
\centering
	\includegraphics[trim=0 0.0cm 0 0.0cm,clip,width=0.7\textwidth]{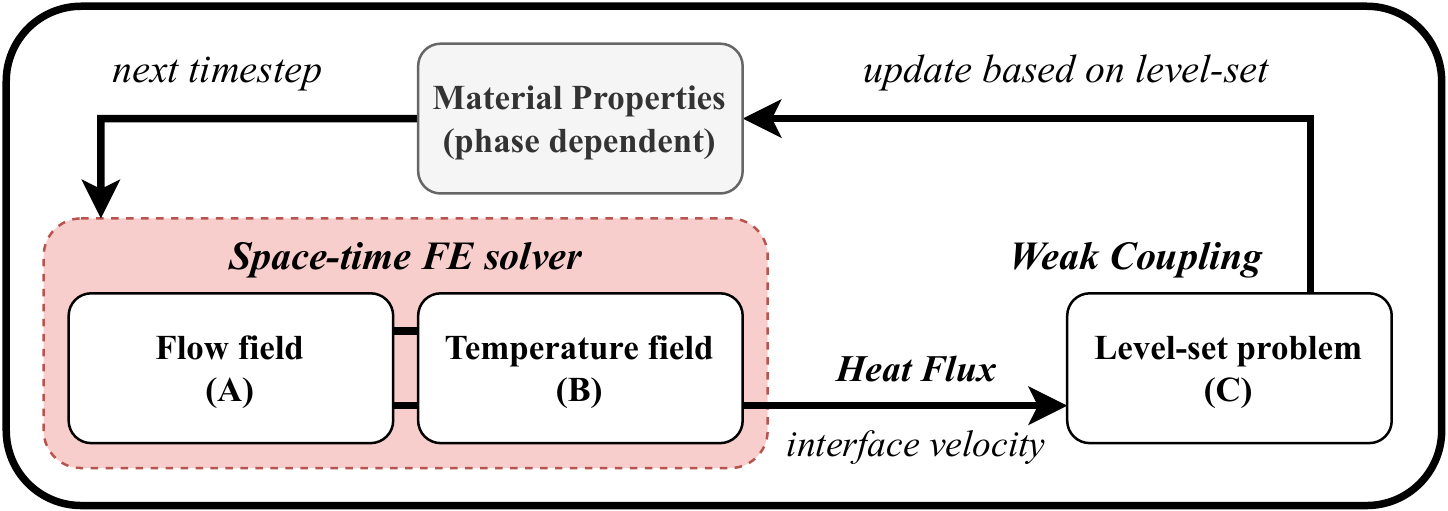}
	\caption{Overview of the numerical approach. Flow field evolution (A) and temperature evolution (B) are computed via space-time finite elements, based on knowledge of the current phase at each point of the domain. The temperature field allows to infer on the evolving phase-change interface (C) by means of the propagating level-set function. This in return allows to update phase information and material properties for the next time step. }
\label{fig:2}
\end{figure}

Our major goal is to determine the evolving flow, temperature and phase distribution, hence the spatio-temporal evolution of velocity field $\mathbf{u}(\mathbf{x},t)=(u(\mathbf{x},t),v(\mathbf{x},t))^\intercal$, pressure field $p(\mathbf{x},t)$ and temperature field  $T(\mathbf{x},t)$.\\
Due to the density being constant within each phase, we are facing an incompressible scenario, for which many numerical methods exist. The fundamental challenge, however, is to accurately account for the evolving PCI. This is necessary as the material parameters differ with the phase and flow occurs in the liquid region only. \\
Our fundamental computational approach combines an established space-time FEM to solve for flow and temperature evolution with a level-set technique to account for the evolving PCI into a novel space-time finite element level-set solver for convection-coupled phase-change processes. A general overview of its essential building blocks (\textbf{A}), (\textbf{B}) and (\textbf{C}) is given in Figure \ref{fig:2}. First, flow field (\textbf{A}), and temperature field (\textbf{B}) are solved by means of a phase-wise space-time finite element dicretization, as detailed in Section \ref{sec:modelling}. Based on that, we determine the evolution of the PCI (\textbf{C}) as a propagating level-set function. Its integration into the space-time finite element framework is detailed in Section \ref{section3}. The PCI evolution allows an update of the liquid and solid region along with an update of the corresponding material properties. From now on, the bold capital letters (\textbf{A}), (\textbf{B}), (\textbf{C}) refer to the algorithms building blocks in Fig.\ \ref{fig:2}. The numerical approaches will be described in more details in the Sections \ref{sec:modelling} and \ref{section3}.

%% file: section2.tex
\section{Space-time finite element approximation to the flow and temperature evolution}
\label{sec:modelling}
This section provides details on the numerical solver for the flow field and the temperature field, hence building blocks (\textbf{A}) and (\textbf{B}) in Figure \ref{fig:2}.

\subsection{Flow field modelling} 
\label{subsec:flowfield}

Let $\Omega\subset\mathbb{R}^2$ be the bounded domain in Figure \ref{fig:basicDomain} and let $t\in(0,T)$ be a time instant. For a two-phase problem, we consider the subdomains $\Omega_1(t)$, $\Omega_2(t)$, such that $\Omega_1(t)\cup\Omega_2(t)=\Omega$ for each $t$. We call $\partial\Omega$ the fixed outer boundary, while the PCI is given by $\Gamma_{\textnormal{PCI}}(t) = \Omega_1(t) \cap \Omega_2(t)$. Note that the outer domain and its boundary do not change with time, while the regions associated with each phase are time-dependent. To compute the flow and pressure fields, see block (\textbf{A}) in Fig.\ \ref{fig:2}, we consider the incompressible Navier-Stokes equations

\begin{alignat}{2}
	\rho_*\left(\frac{\partial\textbf{u}}{\partial t} + \textbf{u}\cdot\nabla\textbf{u}-\textbf{f}\right) - \nabla\cdot\boldsymbol{\sigma_*} &= 0 \hspace{5mm} &&\text{in} \hspace{2mm} \Omega\times(0,T), \label{eq:navierStokes1} \\
	\nabla\cdot\textbf{u} &= 0  &&\text{in} \hspace{2mm} \Omega\times(0,T). \label{eq:navierStokes2}
\end{alignat}
To increase readability we omit the explicit dependency of each component on space and time $(\textbf{x},t)$. We consider a Newtonian fluid, so that we can write the stress tensor $\boldsymbol{\sigma}$ in Eq.\ \eqref{eq:navierStokes1} as
\begin{equation}
	\boldsymbol{\sigma}_*(\textbf{u},p)=-p\mathbf{I} + 2\mu_*\boldsymbol{\varepsilon}(\mathbf{u}),
\end{equation}
where
\begin{equation}
	\boldsymbol{\varepsilon}(\mathbf{u}) = \frac{1}{2}\big(\nabla\textbf{u}+(\nabla\textbf{u})^\intercal\big).
\end{equation}
The subscript $*$, present in the density $\rho_*$ and the dynamic viscosity $\mu_*$, indicates the phase-dependent material properties associated with each subdomain $\Omega_i(t)$, such that
\begin{equation}
	\rho_*(\textbf{x},t) = \begin{cases}
		\rho_1, \hspace{2mm} \text{if} \hspace{2mm} \textbf{x}\in\Omega_1(t), \\
		\rho_2, \hspace{2mm} \text{if} \hspace{2mm} \textbf{x}\in\Omega_2(t),
	\end{cases}
	\mu_*(\textbf{x},t) = \begin{cases}
		\mu_1, \hspace{2mm} \text{if} \hspace{2mm} \textbf{x}\in\Omega_1(t), \\
		\mu_2, \hspace{2mm} \text{if} \hspace{2mm} \textbf{x}\in\Omega_2(t).
	\end{cases}
\label{eq:properties}
\end{equation}
Recall that material properties are constant within each phase. To close the problem we assign Dirichlet and Neumann type boundary conditions and the initial condition
\begin{equation}
	\begin{aligned}
		\textbf{u}(\textbf{x},0) = \textbf{u}_0(\textbf{x}) \hspace{5mm} &\textnormal{in}\hspace{2mm}\Omega,\\
		\textbf{u}=\textbf{g}  \hspace{5mm} &\textnormal{on}\hspace{2mm}\Gamma_g, \\
		\boldsymbol{\sigma}_*\cdot\textbf{n}=\textbf{h}  \hspace{5mm} &\textnormal{on}\hspace{2mm}\Gamma_h.
	\end{aligned}
\label{eq:boundaryConditions}
\end{equation}
The two boundaries $\Gamma_g$ and $\Gamma_h$ denote the parts where we assign Dirichlet and Neumann conditions respectively, such that $\mathring{\Gamma}_g\cap\mathring{\Gamma}_h=\emptyset$ and $\Gamma_g\cup\Gamma_h=\partial\Omega$. We refer to \cite{Salsa2016} for the weak formulation of the problem.

\begin{figure}
    \centering
    \includegraphics[trim={12cm 11cm 12cm 13cm},width=0.5\textwidth]{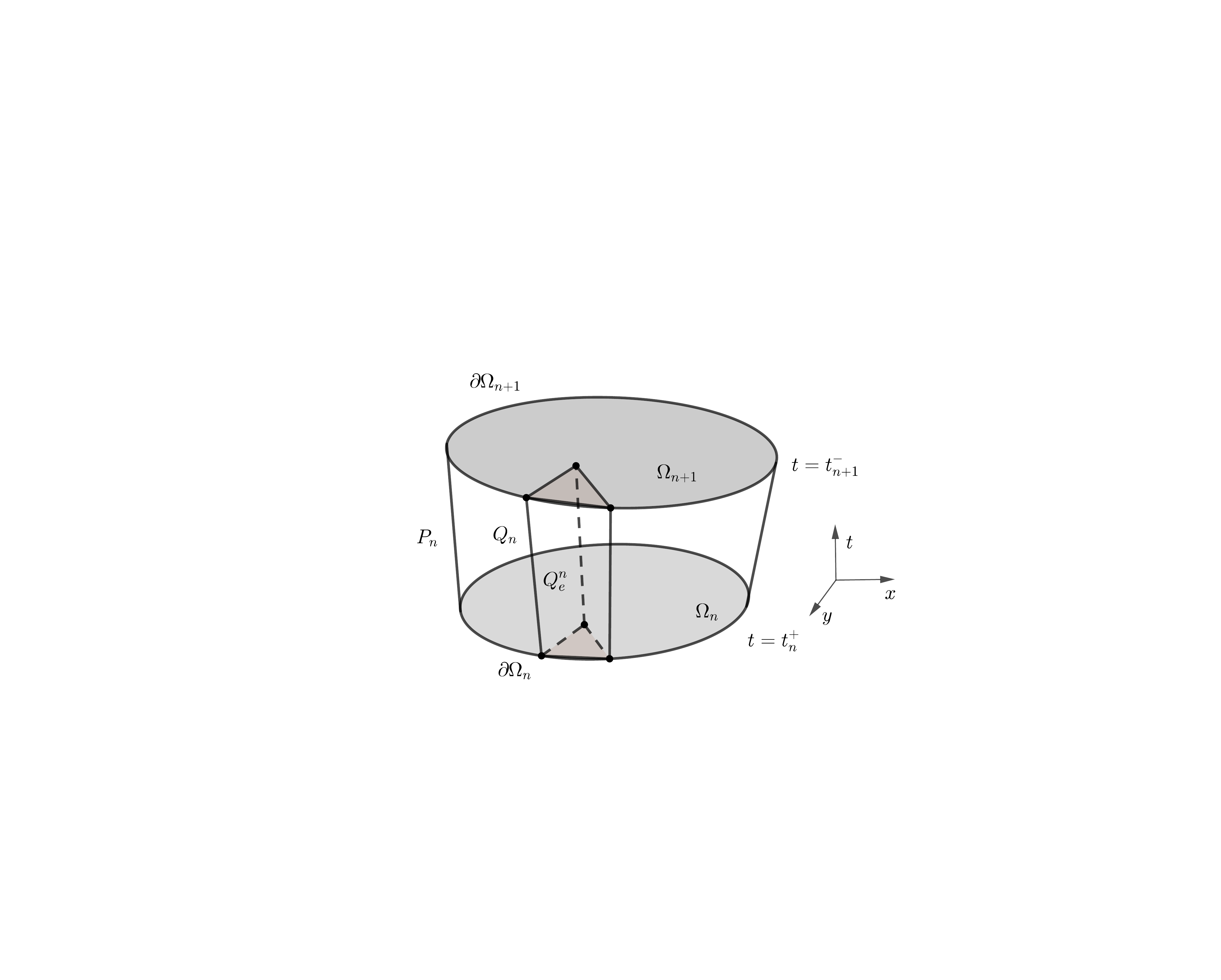}
    \caption{Sketch of a space-time slab $Q_n$ between the lower time level $t_n^+$, with spatial domain $\Omega_n$, and the upper time level $t_{n+1}^-$, with spatial domain $\Omega_{n+1}$. We call the lateral surface $P_n$ and the generic space-time element $Q_n^e$.}
    \label{fig:spaceTimeDomain}
\end{figure}

To solve the time-dependent Navier-Stokes system we employ the space-time FEM. Instead of considering the weak form of the equations only over the spatial domain, we use a specific domain in space and time \cite{Huerta2005}. Let us divide the time interval $(0,T)$ into subintervals $I_n=(t_n,t_{n+1})$, such that $0=t_0<t_1<\dots<t_N=T$. At the time level $t_n$ we have the spatial domain $\Omega_n=\Omega_{t_n}$ and its boundary $\partial\Omega_n=\partial\Omega_{t_n}$. Then, we define the space-time slab $Q_n$ as the domain enclosed by the surfaces $\Omega_n,\Omega_{n+1}$ and the surface $P_n$. The latter is described by $\partial\Omega_t$ when $t$ passes $I_n$, see Figure \ref{fig:spaceTimeDomain}. For each space-time slab, we define the interpolation and weighting function spaces for first order polynomials
\begin{equation}
\begin{aligned}
	(\mathcal{S}_{\textbf{u}}^h)_n &= \{\textbf{u}^h\in[H^{1h}(Q_n)]^2 \hspace{1mm}|\hspace{1mm} \textbf{u}^h=\textbf{g} \hspace{1mm}\text{on}\hspace{1mm}  (P_n)_g\}, \\
	(\mathcal{V}_{\textbf{u}}^h)_n &= \{\textbf{w}^h\in[H^{1h}(Q_n)]^2 \hspace{1mm}|\hspace{1mm} \textbf{w}^h=\textbf{0} \hspace{1mm}\text{on}\hspace{1mm}  (P_n)_g\}, \\
	(\mathcal{S}_p^h)_n &= (\mathcal{V}_p^h)_n=\{p^h\in H^{1h}(Q_n)\},
\end{aligned}
\label{eq:NSfunctionSpace}
\end{equation}
where the subscripts $\textbf{u},p$ indicate pressure and velocity. $(P_n)_g$ and $(P_n)_h$ denote the portions of the space-time boundary of Dirichlet and Neumann type, respectively. Having defined the functional setting, we follow the formulation presented in \cite{Pauli2017.2}. The stabilized space-time formulation for the incompressible Navier-Stokes Equations \eqref{eq:navierStokes1} and \eqref{eq:navierStokes2} then reads: \\ \textit{Given $(\textbf{u}^h)^-_n$, find $\textbf{u}^h\in(\mathcal{S}_{\textbf{u}}^h)_n$ and $p^h\in(\mathcal{S}_p^h)_n$ such that $\forall\textbf{w}^h\in(\mathcal{V}^h_{\textbf{u}})_n$ and $\forall q^h\in(\mathcal{V}_p^h)_n$}
\begin{multline}
	\int_{Q_n}\textbf{w}^h\cdot\rho_*\left(\frac{\partial\textbf{u}^h}{\partial t} + \textbf{u}^h\cdot\nabla\textbf{u}^h - \textbf{f}\right) \diff Q + \int_{Q_n}\nabla\textbf{w}^h:\boldsymbol{\sigma}_*(\textbf{u}^h,p^h)\diff Q \\ + \int_{Q_n}q^h\nabla\cdot\textbf{u}^h\diff Q + \int_{\Omega_n}(\textbf{w}^h)^+_n\cdot\rho_*\left((\textbf{u}^h)^+_n - (\textbf{u}^h)^-_n\right)\diff\Omega \\ + \sum_{e=1}^{(n_{\text{el}})_n}\int_{Q_n^e}\tau_{\text{MOM}}\frac{1}{\rho_*}\left[\rho_*\left(\textbf{u}^h\cdot\nabla\textbf{w}^h\right) + \nabla q^h\right] \\ \cdot \left[\rho_*\left(\frac{\partial\textbf{u}^h}{\partial t} + \textbf{u}^h\cdot\nabla\textbf{u}^h-\textbf{f}\right) - \nabla\cdot\boldsymbol{\sigma}_*(\textbf{u}^h,p^h)\right]\diff Q \\ + \sum_{e=1}^{(n_{\text{el}})_n}\int_{Q^e_n}\tau_{\text{CONT}}\nabla\cdot\textbf{w}^h\rho_*\nabla\cdot\textbf{u}^h\diff Q = \int_{(P_n)_h}\textbf{w}^h\cdot\textbf{h}^h\diff P. \label{eq:spaceTimeNS}
\end{multline}
Note that in the above equation we employ the notation
\begin{align}
	(\textbf{u}^h)^{\pm}_n & =\lim_{\epsilon\rightarrow 0}\textbf{u}(t_n\pm\epsilon), \nonumber \\
	\int_{Q_n}(\dots)\diff Q &= \int_{I_n}\int_{\Omega_t}(\dots)\diff\Omega\diff t, \label{eq:STnotation}\\
	\int_{P_n}(\dots)\diff P &= \int_{I_n}\int_{\partial\Omega_t}(\dots)\diff\Gamma\diff t. \nonumber
\end{align}
The system is solved sequentially for each space-time slab, starting with $(\textbf{u}^h)^+_0 = \textbf{u}_0$. The fourth term in Eq.\ \eqref{eq:spaceTimeNS} is the jump term, which induces the weak continuity in time for the velocity field over $\Omega_n$ \cite{Huerta2005}. The two terms $\tau_{\text{MOM}}$ and $\tau_{\text{CONT}}$ stabilize the momentum and continuity equations, respectively. The expressions for their values can be found in \cite{Pauli2017.2}.

\subsection{Temperature evolution}
\label{subsec:tempfield}

To model the temperature field $T(\textbf{x},t)$, see block (\textbf{B}) in Fig.\ \ref{fig:2}, we consider the transient heat equation
\begin{equation}
	\rho_*(c_p)_*\left(\frac{\partial T}{\partial t} + \textbf{u}\cdot\nabla T\right)=\kappa_*\Delta T \hspace{5mm} \text{in} \hspace{2mm} \Omega\times(0,T), \label{eq:heat}
\end{equation}
where $(c_p)_*$ is the heat capacity and $\kappa_*$ is the thermal conductivity, that again might vary with the phase. The term $\textbf{u}(\textbf{x},t)$ is the velocity field from Eqs.\ \eqref{eq:navierStokes1}, \eqref{eq:navierStokes2}. We consider again Dirichlet and Neumann type boundary conditions to close the problem
\begin{equation}
	\begin{aligned}
		T(\textbf{x},0) = T_0(\textbf{x}) \hspace{5mm} &\textnormal{in}\hspace{2mm}\Omega,\\
		T=g  \hspace{5mm} &\textnormal{on}\hspace{2mm}\Gamma_g, \\
		\nabla T\cdot\textbf{n}=h  \hspace{5mm} &\textnormal{on}\hspace{2mm}\Gamma_h.
	\end{aligned}
\end{equation}
Notice that the boundaries $\Gamma_g$ and $\Gamma_h$ in the temperature equation can differ from $\Gamma_g$ and $\Gamma_h$ in the Navier-Stokes equations. An in depth analysis of the weak formulation of Problem \eqref{eq:heat} can be found in \cite{Salsa2016}. In order to introduce the space-time formulation, we modify the functional space in Eq.\ \eqref{eq:NSfunctionSpace} for a scalar problem, such that
\begin{equation}
\begin{aligned}
	(\mathcal{S}_T^h)_n &= \{T^h\in H^{1h}(Q_n) \hspace{1mm}|\hspace{1mm} T^h=g \hspace{1mm}\text{on}\hspace{1mm}  (P_n)_g\}, \\
	(\mathcal{V}_T^h)_n &= \{v^h\in H^{1h}(Q_n) \hspace{1mm}|\hspace{1mm} v^h=0 \hspace{1mm}\text{on}\hspace{1mm}  (P_n)_g\}. \\
\end{aligned}
\label{eq:heatFunctionSpace}
\end{equation}
The stabilized space-time formulation of Eq.\ \eqref{eq:heat} reads: \\ \textit{Given $(T^h)^-_n$, find $T^h\in(\mathcal{S}_T^h)_n$ such that $\forall v^h\in(\mathcal{V}^h_T)_n$}
\begin{multline}
	\int_{Q_n}v^h\cdot\rho_*(c_p)_*\left(\frac{\partial T^h}{\partial t} + \textbf{u}^h\cdot\nabla T^h\right)\diff Q + \int_{Q_n} \nabla v^h\cdot\kappa_*\nabla T^h\diff Q \\ + \int_{\Omega_n}(v^h)^+_n\rho_*(c_p)_*[(T^h)^+_n - (T^h)^-_n]\diff\Omega \\
	= \sum_{e=1}^{(n_{\text{el}})_n} \int_{Q_n^e} \tau_{\text{TEMP}}\frac{1}{\rho_* (c_p)_*} \left[\rho_* (c_p)_* \left( \frac{\partial v^h}{\partial t} + \textbf{u}^h \cdot\nabla v^h\right)\right] \\\cdot \left[\rho_*(c_p)_* \left(\frac{\partial T^h}{\partial t} + \textbf{u}^h\cdot\nabla T^h\right) -\nabla\cdot\kappa_*\nabla T^h\right] \diff Q = \int_{(P_n)_h} v^hh^h\diff P,
\label{eq:heatSpaceTime}
\end{multline}
with $(T^h)^+_0=T_0$. The notation is analogous to the one introduced with Eq.\ \eqref{eq:STnotation} and we refer to \cite{Pauli2017.2} for the value of the stabilization term $\tau_{\text{TEMP}}$. Note that problem (\textbf{B}) is coupled with problem (\textbf{A}) through the advection velocity in Eq.\ \eqref{eq:heat}. This constitutes a one-way coupling as  we do not consider temperature contributions in the Navier-Stokes Eqs.\ \eqref{eq:navierStokes1}, \eqref{eq:navierStokes2}, such as a bouyancy term.

%% file: section3.tex
\section{Level-set approach to tracking the phase-change interface}
\label{section3}

\begin{figure}

    \begin{subfigure}[c]{0.49\textwidth}
    \centering
    \includegraphics[trim={0.0cm 0.0cm 0.0cm 0.0cm},clip, width=0.7\textwidth]{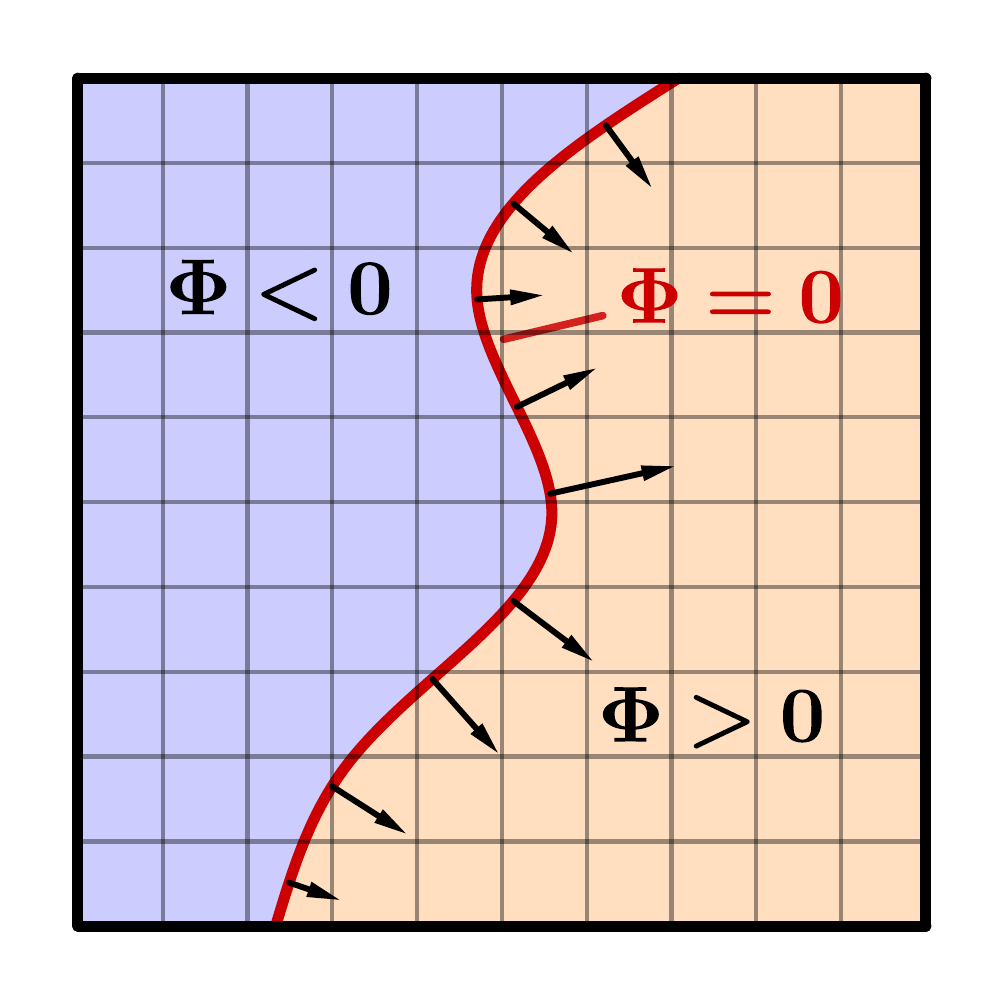}
        \subcaption{$t=t_1>0$}
        \label{fig:levvel}
    \end{subfigure}
    \begin{subfigure}[c]{0.49\textwidth}
    \centering
    \includegraphics[trim={0.0cm 0.0cm 0.0cm 0.0cm},clip, width=0.7\textwidth]{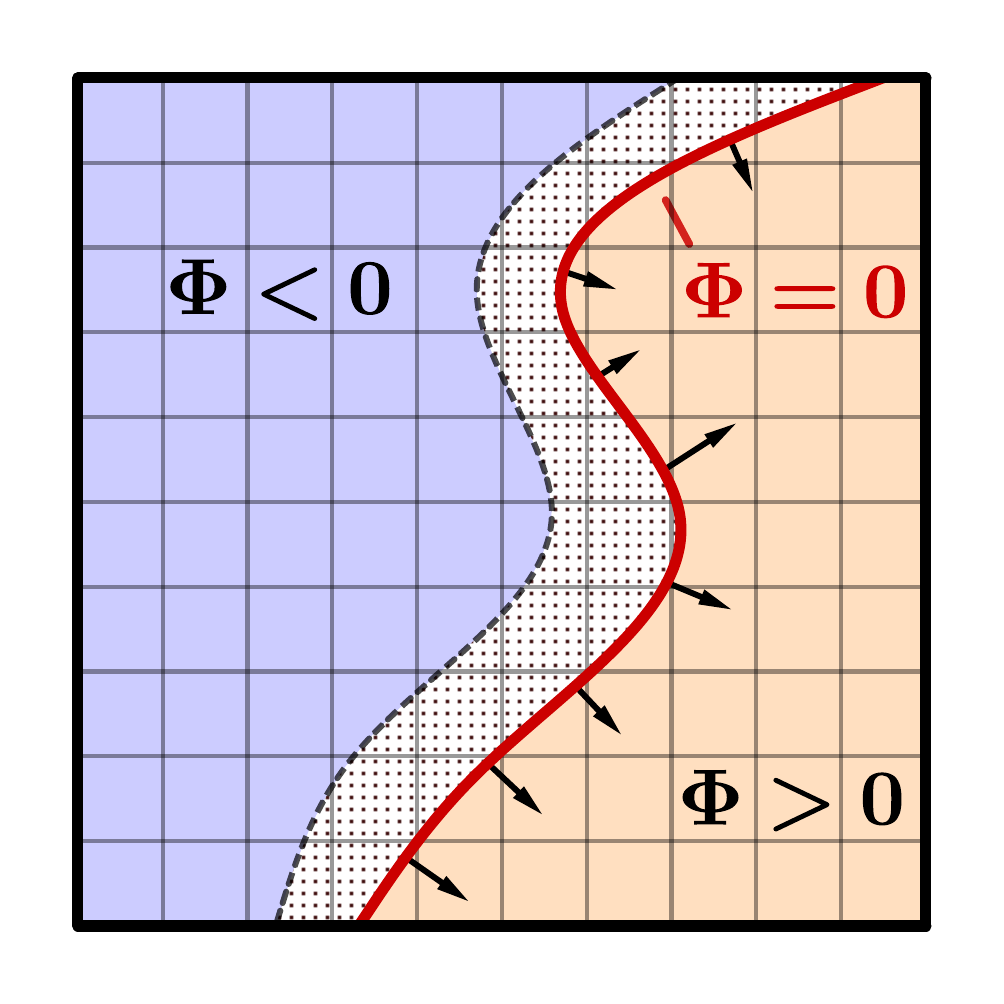}
        \subcaption{$t=t_2 > t_1$}
        \label{fig:levshift}
    \end{subfigure}

    \caption{Evolution of the phase-change interface (PCI), depicted in red. Due to melting or solidification, the subdomains associated with each phase might change over time. On the left, the computational domain is shown at time $t_1$. On the right, the updated domain is shown at time $t_2>t_1$. Note that the PCI shifts to the right and the left phase ($\phi < 0$) acquires the dotted region. The black arrows represent the local propagation speed of the interface, that is the advective term $\textbf{v}(\textbf{x},t)$ in Eq.\ \eqref{eq:levelSet}. Such velocity is modelled as the Stefan condition.}
    
    \label{fig:levelSetVelocity}
    
\end{figure}

In this section we describe the level-set method to handle the evolving phase-change interface, that is building block (\textbf{C}) in Figure \ref{fig:2}. Then, we introduce the Stefan condition to model the local propagation velocity of the interface. Finally, we discuss the reconstruction of the heat-flux discontinuity at the interface via the ghost-split approach.

\subsection{General formulation of the level-set method}
Let $\Phi:\Omega\times(0,T)\rightarrow\mathbb{R}$ be a scalar function. The function $\Phi$ is continuous and it is defined as
\begin{equation}
	\begin{aligned}
		\Phi(\textbf{x},t)<0 \hspace{5mm} &\textnormal{in}\hspace{2mm}\Omega_1(t),\\
		\Phi(\textbf{x},t)>0 \hspace{5mm} &\textnormal{in}\hspace{2mm}\Omega_2(t),\\
		\Phi(\textbf{x},t)=0  \hspace{5mm} &\textnormal{on}\hspace{2mm}\Gamma_{\text{PCI}}(t).
	\end{aligned}
\end{equation}
This function is called the level-set function, because the interface $\Gamma_{\text{PCI}}(t)$ is its zero level set, that is
\begin{equation}
	\Gamma_{\text{PCI}}(t) = \{\textbf{x}\in\Omega : \Phi(\textbf{x},t)=0\}.
\end{equation}
As such, the function $\Phi$ indicates in which subdomain a point $\textbf{x}$ is located. Then, the material properties can be expressed as function of $\Phi$. For instance the density $\rho_*$ in Eq.\ \eqref{eq:properties} can be written as
\begin{equation}
	\rho_* = \rho_1 + (\rho_2 - \rho_1)H_\epsilon(\Phi),
\end{equation}
so that Eqs.\ \eqref{eq:navierStokes1}, \eqref{eq:navierStokes2} and \eqref{eq:heat} can describe two different phases. To avoid sharp changes in the material properties across the PCI, we select the function $H_\epsilon(\cdot)$ as the smoothed Heaviside function
\begin{equation}
	H_\epsilon(\Phi) = \begin{cases} 0, &\Phi < -\epsilon, \\ 
	\frac{1}{2}\left[1 + \frac{\Phi}{\epsilon} + \frac{1}{\pi}\sin\left(\frac{\pi\Phi}{\epsilon}\right)\right], \hspace{3mm} &|\Phi| \leq \epsilon, \\
	1, &\Phi > \epsilon, \end{cases}
	\label{eq:heaviside}
\end{equation}
for some fixed, small $\epsilon$. The interface now has a fixed thickness of approximately $\frac{2\epsilon}{|\nabla\Phi|}$, which is proportional to the spatial mesh size. We refer to \cite{sethian2003} for a detailed discussion of this formulation and its advantages. The evolution of $\Gamma_{\text{PCI}}(t)$ is described by the equations
\begin{equation}
	\begin{aligned}
		\frac{\partial\Phi}{\partial t} + \textbf{v}\cdot\nabla\Phi = 0 \hspace{5mm}&\textnormal{in}\hspace{2mm}\Omega\times(0,T), \\
		\Phi(\textbf{x},0) = \Phi_0(\textbf{x}) \hspace{5mm}&\textnormal{in}\hspace{2mm}\Omega,
	\end{aligned}
	\label{eq:levelSet}
\end{equation}
where $\textbf{v}$ denotes the propagation velocity of the PCI. The initial condition $\Phi_0(\textbf{x})$ is chosen such that $\Phi(\textbf{x},t)$ is the signed distance function with respect to the PCI. Note that we have obtained a scalar advection problem, which shares many similarities with Eq.\ \eqref{eq:heat}. Thus, we do not repeat the details of the space-time formulation and the solution approach. Additional information can be found in Section 3.10 of \cite{Huerta2005}. Also note that the transport term $\textbf{v}$ is not known. This is the major challenge of the model, which we will address in the next sections.

By construction, the interface stays sharp and the determination of its normals and its curvature are straightforward \cite{Quarteroni2010}. In particular, we have 
\begin{equation}
	\textbf{n}_{\Gamma_{\text{PCI}}} = \frac{\nabla\Phi}{|\nabla\Phi|}, \hspace{5mm} \kappa_{\Gamma_{\text{PCI}}} = -\nabla\cdot\left(\frac{\nabla\Phi}{|\nabla\Phi|}\right),
\label{eq:levelSetProperties}
\end{equation}
where $\textbf{n}_{\Gamma_{\text{PCI}}}$ is the interface unit normal from $\Omega_1(t)$ to $\Omega_2(t)$ and $\kappa_{\Gamma_{\text{PCI}}}$ is the curvature. While solving Eq.\ \eqref{eq:levelSet}, if the gradient of $\Phi$ becomes too large with respect to the grid spacing, we lose accuracy in the interface representation. To avoid it, we reinitialize $\Phi$ using the signed distance with respect to the current PCI after a certain number of time steps \cite{Quarteroni2010}. This entails computing the shortest distance to the interface of all nodal points on the numerical grid. The most naive implementation of this procedure, which is also used in this paper, has a complexity of
\begin{equation}
    \mathcal{O}(n_{\text{grid}} \cdot n_{\text{PCI}}),
\end{equation}
where $n_{\text{grid}}$ is the number of mesh nodes and $n_{\text{PCI}}$ is the number of interface crossings with the mesh. A popular algorithmic alternative that will grant a performance gain in the future would be the fast marching method \cite{Kimmel1998}. Note that the overhead introduced by the reinitialization procedure is negligible compared to the cost of solving for the flow and temperature fields. For the 2D cases considered in this paper, this does not result in a critical performance bottleneck, but it might become relevant for complex 3D meshes.

\subsection{Interface propagation and Stefan condition}

So far, we can compute the evolving flow, pressure and temperature fields and we can track the PCI propagation, yet we do not have a closure for the level-set's advection term $\textbf{v}(\textbf{x},t)$ in Eq.\ \eqref{eq:levelSet}. This effectively corresponds to the need for formulating the coupling of block (\textbf{B}) into (\textbf{C}) in Fig.\ \ref{fig:2}. At the interface itself the propagation velocity $\textbf{v}(\textbf{x},t)$ has to correspond to the actual phase-change rate, hence melting or solidification rate.
This phase-change rate can readily be determined from local energy conservation across the PCI, which gives rise to a heat-flux jump condition also known as the Stefan condition, that is
\begin{equation}
	\rho \, h_m\textbf{U}(\textbf{x},t) = -\kappa_L\nabla T\big|_{\textbf{X}^-} + \kappa_S\nabla T\big|_{\textbf{X}^+} = \left[\kappa\nabla T\right]^S_L = q_L - q_S.
\label{eq:meltingVelocity}
\end{equation}
Here, $h_m$ is the latent heat of melting, $\rho$ denotes the material's density, $\kappa$ the material's conductivity, $\textbf{X}^\pm$ denote the limits taken from either side of the PCI and $[\cdot]^S_L$ refers to the liquid and solid regions. $\textbf{U}(\textbf{x},t)$ finally stands for the Stefan velocity, hence the rate at which the interface changes its phase. The Stefan velocity, which is proportional to the interfacial heat flux, then defines the coupling between the temperature and the level-set equations and provides a closure for the level-set's propagation term $\textbf{v}(\textbf{x},t)\equiv\textbf{U}(\textbf{x},t)$ in Eq.\ \eqref{eq:levelSet}. Its evaluation requires an accurate recovery of the temperature gradient within our space-time FEM framework.

\subsection{Heat-flux reconstruction at the interface} 
\label{heatFluxReconstruction}

\begin{figure}
\centering
    \includegraphics[trim={0.0cm 0.0cm 0.0cm 0.0cm},clip, width=0.7\textwidth]{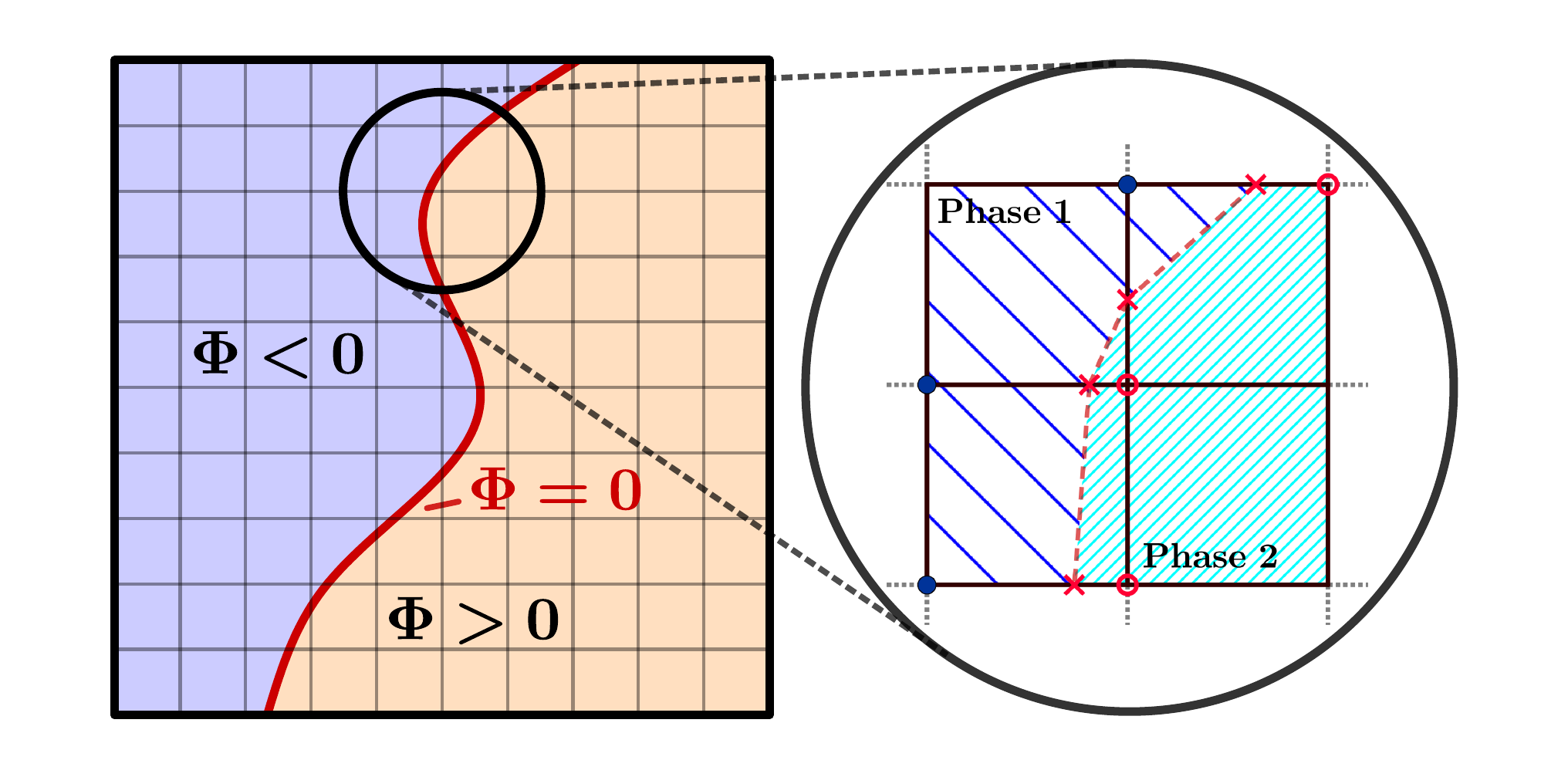}
    \caption{Uniform structured mesh bisected by the PCI. The close-up view shows the intersections with four element faces that are represented by red crosses. The nodes adjacent to the intersection points are marked by blue dots (\textit{Phase 1}) and red circles (\textit{Phase 2}). The numerical PCI results from piecewise linear shape functions.} 
    \label{fig:zoomInterfaceCut}
\end{figure}

In the previous sections we have described all the individual modules of our problem as well as how to couple them. What we are missing is the computation of the jump term in Eq.\ \eqref{eq:meltingVelocity}. We now propose a method to recover the heat-flux discontinuity at the PCI using FEM with element-wise continuously differentiable shape functions. Existing algorithms that deal with jumps in the first derivatives have been described in \cite{Gibou2002, Gibou2005, Gibou2007}, but these publications only consider finite difference discretizations. By extending the idea to our space-time formulation, we can take advantage of the versatility of FEM without the need of local enrichments of the FEM basis. In particular, the total number of degrees of freedom stays constant over time, meaning that we do not alter the size of the global system matrix. This can be an advantage in view of highly parallelized FEM codes.

\begin{figure}

    \centering
  \begin{subfigure}[c]{0.3\textwidth}
    \includegraphics[trim={0.5cm 0.7cm 0.5cm 0.7cm},clip, width=\textwidth]{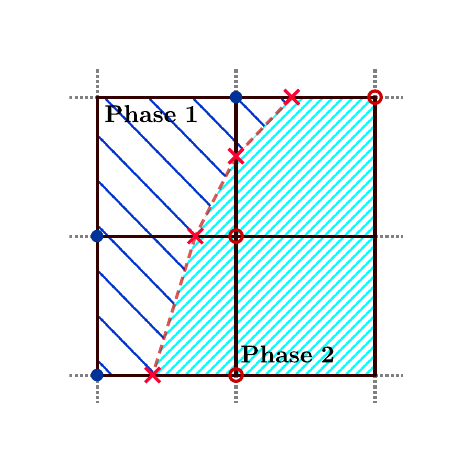}
        \subcaption{}
        \label{fig:fluxnodes1}
    \end{subfigure}
    \begin{subfigure}[c]{0.3\textwidth}
    \includegraphics[trim={0.5cm 0.7cm 0.5cm 0.7cm},clip, width=\textwidth]{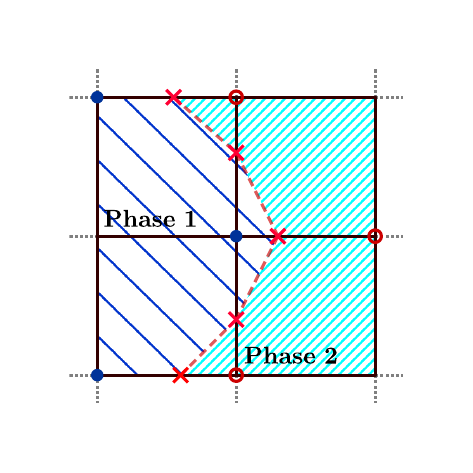}
        \subcaption{}
        \label{fig:fluxnodes2}
    \end{subfigure}
    \begin{subfigure}[c]{0.3\textwidth}
    \includegraphics[trim={0.5cm 0.7cm 0.5cm 0.7cm},clip, width=\textwidth]{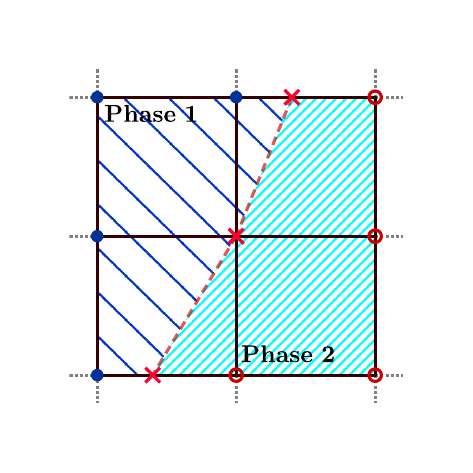}
        \subcaption{}
        \label{fig:fluxnodes4}
    \end{subfigure}
    
    \caption{Recovery of the flux discontinuity on a 2D structured mesh. The red crosses mark the intersection points of the interface (red dashed line) with element faces. The highlighted element nodes are involved in the computation of the heat flux at the intersection point cutting an adjacent face. The blue and red circles indicate nodes belonging to phase 1 (left) and to phase 2 (right), respectively. We show examples for various shapes of the interface. Figure \ref{fig:fluxnodes4} shows a special case where a mesh node is located on the PCI.}
    
    \label{fig:fluxNodes}
    
\end{figure}

For the sake of a simpler visualization we describe the method on a uniform structured mesh, but note that all the concepts are applicable to the unstructured case. In fact, multiple examples on unstructured grids are provided in Section \ref{section4}. Figure \ref{fig:zoomInterfaceCut} gives a close-up view of some elements bisected by the PCI in our discretized domain. The intersection points, shown as red crosses, are described by the equation $\Phi(\textbf{x},t)=0$. The first issue comes from the choice of evaluation points for the representative fluxes $q_L$, $q_S$ in Eq.\ \eqref{eq:meltingVelocity}. Let us consider various shapes of the interface as shown in Figure \ref{fig:fluxNodes}. The key requirement is that the flux converges to the value at the interface in the limit for fine mesh resolutions. The evaluation of the temperature at an arbitrary point in the domain is not straightforward, since it requires a mapping from physical coordinates to local element coordinates. Similarly, the evaluation at points normal to the interface presents issues, since the normal is not well-defined at intersection points (red crosses in Fig.\ \ref{fig:fluxNodes}). To circumvent these difficulties we adopt a strategy based on three propositions:

\begin{enumerate}
	\item If an element face is cut by the PCI, we consider the nodes that belong to this face as \textit{flux nodes}, which means that we use the numerical gradient at these nodes as the representative fluxes in Eq.\ \eqref{eq:meltingVelocity} to compute the Stefan velocity at the crossing. Based on the nodal value of the level-set function we can determine the associated phase of a node, yielding the flux $q_L$ or $q_S$;
	\item We compute each nodal gradient with a least-squares fitting based on the gradient within all elements adjacent to a mesh node. For example, the central node in Figure \ref{fig:fluxNodes} gives an average of the gradient at the four square elements surrounding it. In case of piecewise linear interpolation functions, this reduces to averaging the element-wise constant gradients;
	\item If the PCI intersects a mesh node, we consider the average of all adjacent nodes in each face to obtain $q_L$ and $q_S$, see Figure \ref{fig:fluxnodes4}.
\end{enumerate}

\begin{figure}
\centering
    \includegraphics[trim={0.0cm 0.0cm 0.0cm 0.0cm},clip, width=0.6\textwidth]{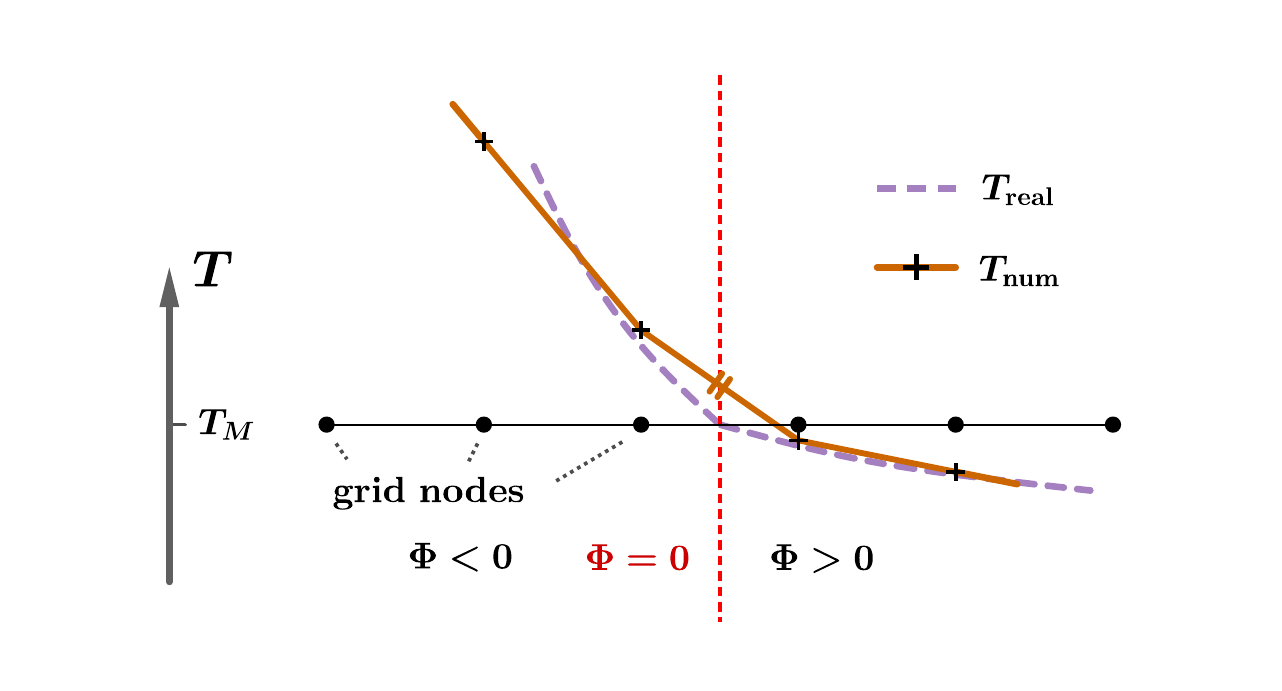}
    \caption{Piecewise linear approximation $T^h(\textbf{x},t)$ (orange) of a fictitious temperature profile $T_{\text{real}}$ (dashed purple), $T_M$ is the phase-change temperature. The central element of the grid is bisected by the PCI, with grid nodes located to the left ($\Phi(x) < 0$) and right ($\Phi(x) > 0$) belonging to different phases. The numerical approximation is element-wise differentiable and does not represent the discontinuous heat flux at the interface.}
    \label{fig:temp1Dnum}
\end{figure}

In order to evaluate the temperature gradient and close the problem, we use the numerical approximation $T^h(\textbf{x},t)$ of the temperature field that we obtain from Eq.\ \eqref{eq:heatSpaceTime}. Following our finite element formulation, the numerical solution can be written as

\begin{equation}
	T^h(\textbf{x},t)=\sum_{i=1}^{n_{\text{grid}}} N_i(\textbf{x}) \, u^n_{i,T},
\end{equation}
where $n_{\text{grid}}$ is the number of nodes, $N_i(\textbf{x})$ and $u^n_{i,T}$ represent the FEM interpolation functions and the nodal temperature values, respectively. Then we shift the derivative onto the interpolation functions such that
\begin{equation}
	\nabla T^h(\textbf{x},t) = \sum_{i=1}^{n_{\text{grid}}} \nabla N_i(\textbf{x}) \, u_{i,T}^n.
\end{equation}
It is clear that the mathematical properties of the approximate gradient $\nabla T^h(\textbf{x},t)$ depend on the properties of the space $\mathcal{S}^h$. In particular, one needs to observe that discontinuities in the gradient can in general not be considered, as the employed piecewise linear interpolation functions are continuous within the elements. Figure \ref{fig:temp1Dnum} shows a fictitious 1D temperature profile where the exact solution, depicted by the dashed purple line, features a discontinuity at the PCI. Such discontinuity is not captured by the numerical solution, shown in orange, which would give $U(x,t)=(\rho h_m)^{-1} \partial_x T|_L^S=0$ across the element cut by the PCI. We will address this matter in Section \ref{ghostSplit} with the ghost-split method.

\subsection{The ghost-split method}
\label{ghostSplit}

\begin{figure}
\centering
    \begin{subfigure}[c]{0.40\textwidth}
    \includegraphics[trim={0.0cm 0.0cm 1.0cm 0.0cm},clip, width=\textwidth]{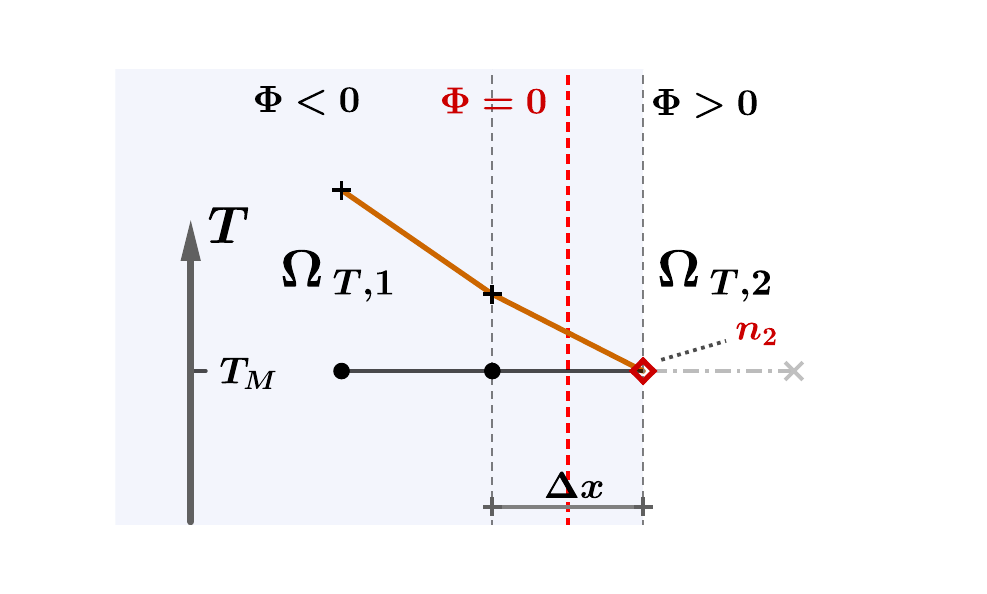}
        \subcaption{Ghost-split domain for the left phase}
        \label{fig:ghost1dleft}
    \end{subfigure}
    \begin{subfigure}[c]{0.40\textwidth}
    \includegraphics[trim={0.0cm 0.0cm 1.0cm 0.0cm},clip, width=\textwidth]{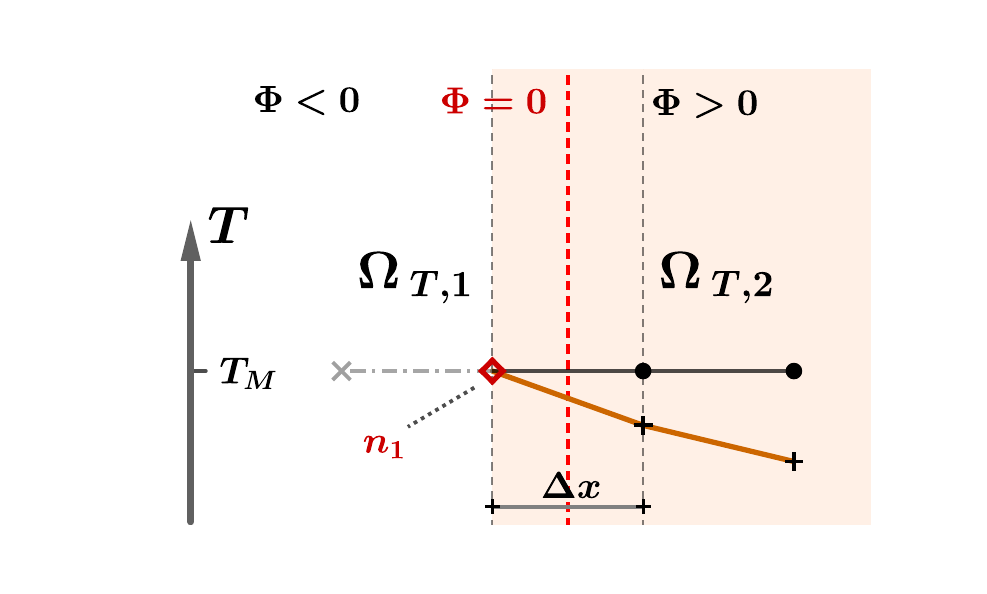}
        \subcaption{Ghost-split domain for the right phase}
        \label{fig:ghost1dright}
    \end{subfigure}

     \caption{1D ghost split of Figure \ref{fig:temp1Dnum}. We split the domain into two subdomains $\Omega_{T,1}$ and $\Omega_{T,2}$ based on the phase information from the level-set function. The two temperature problems are independently solved, but they share the Dirichlet \textit{ghost nodes} (red rhombi), where we prescribe the melting temperature $T_m$. In this example $n_2$ is the ghost node for the domain $\Omega_{T,1}$ (\ref{fig:ghost1dleft}) and $n_1$ is the ghost node for $\Omega_{T,2}$ (\ref{fig:ghost1dright}).}
    \label{fig:ghost1d}
    
\end{figure}

\begin{figure}

\centering
    \begin{subfigure}[c]{0.40\textwidth}
    \includegraphics[trim={0.0cm 0.0cm 1.0cm 0.0cm},clip, width=\textwidth]{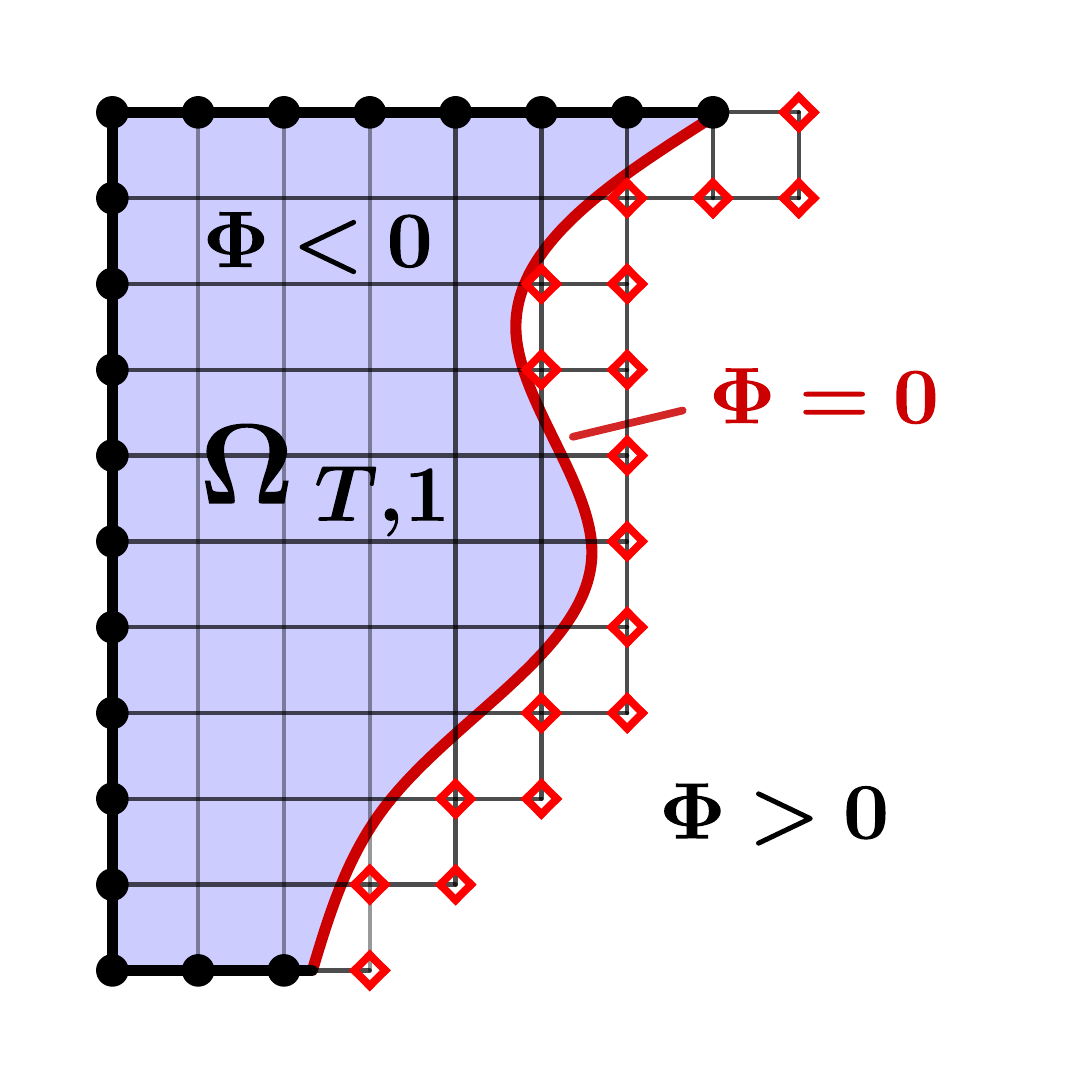}
        \subcaption{Subdomain of the first phase, $\Phi<0$}
        \label{fig:ghost2dleft}
    \end{subfigure}
    \begin{subfigure}[c]{0.40\textwidth}
    \includegraphics[trim={0.0cm 0.0cm 1.0cm 0.0cm},clip, width=\textwidth]{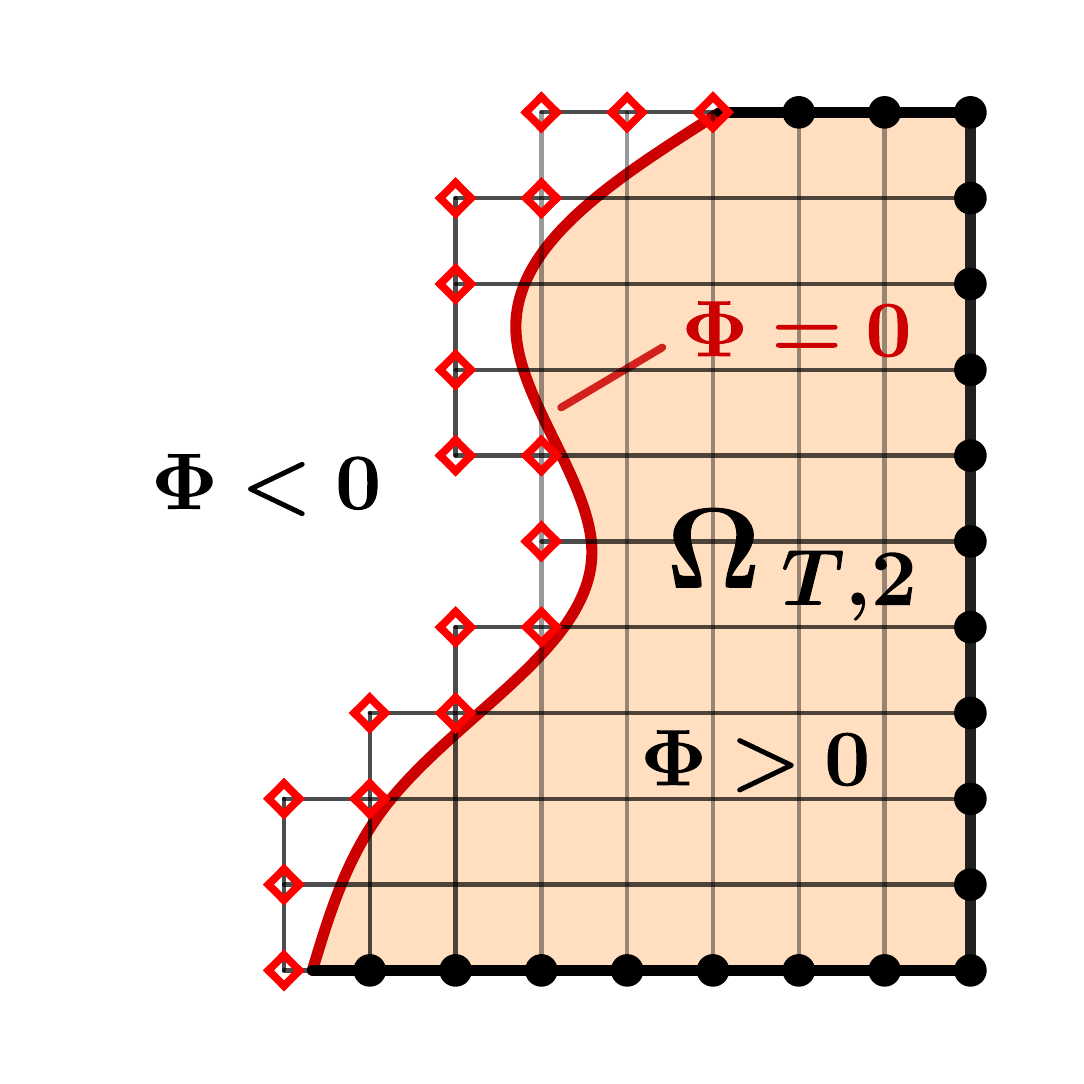}
        \subcaption{Subdomain of the second phase, $\Phi>0$}
        \label{fig:ghost2dright}
    \end{subfigure}

     \caption{2D ghost split of the computational domain. The domain $\Omega$ is divided into subdomains $\Omega_{T,1}$ and $\Omega_{T,2}$ based on the value of the level-set function $\Phi(\textbf{x},t)$. The melting temperature $T_m$ is imposed as a Dirichlet condition for the two subproblems at the \textit{ghost nodes} shown as red rhombi.}
    \label{fig:ghost2d}
    
\end{figure}

Our choice of interpolation functions is not able to represent discontinuous first derivatives in elements cut by the PCI. As discussed in Section \ref{heatFluxReconstruction}, piecewise linear shape functions can show discontinuities in the temperature gradient only at element nodes. The ghost-split method is founded upon the idea that the temperature field associated to each phase can be treated independently. Since the process of melting (or solidification) requires that the numerical solution equals the melting temperature $T_m$ at the PCI, i.e.\ $T^h(\textbf{x}_{\text{PCI}},t) = T_m$ at all times $t$, the interface can be viewed as a Dirichlet type boundary for the adjacent phases. Given this condition, we solve the heat equation in each subdomain $\Omega_i(t)$ without knowledge of the temperature profile in the other phase. Then, we use the gradients recovered form the subproblems to compute the interface propagation velocity as shown in Section \ref{heatFluxReconstruction}. However we can only impose boundary conditions on mesh nodes, which gives rise to the concept of \textit{ghost nodes}. The term refers to the fact that we add additional nodes to the subdomain to enforce the melting temperature at the approximate position of the interface. Figure \ref{fig:ghost1d} shows an example for our fictitious 1D case. Note that an additional node, depicted as a red rhombus, is added to each phase and the melting temperature $T_m$ is imposed on it. Figure \ref{fig:ghost2d} extends the method to our introductory 2D domain of a two-phase problem. Even if we limit ourselves to the description of a 2D case, the algorithm does not depend on the number of space dimensions. Note that the overall number of nodes on which we solve for the temperature field does not change compared to the original discretization shown in Fig.\ \ref{fig:zoomInterfaceCut}.

The ghost-split technique enforces the melting temperature at a node close to the PCI, which introduces an error in the interface location computed at later time steps. Following the argument in \cite{Gibou2002} we note that this error corresponds to a $\mathcal{O}(\Delta x)$ disturbance in the interface position. In other words, imposing the melting temperature at a ghost node misplaces the PCI by at most one element length. We observe however that as the mesh gets finer, the position of the ghost node converges to the correct interface location, that is $\Delta x\rightarrow0$. Note that unlike the structured grid used in \cite{Gibou2002}, the element face length may vary across the mesh on arbitrary unstructured grids.

\subsection{Time step control and temperature extrapolation}
\label{temperatureExtrapolation}

Before we show our numerical results, we need to address two more aspects. Let us consider again the 1D case of a two-phase problem as in Fig.\ \ref{fig:ghost1d}. Following our procedure we compute the advection velocity $\textbf{v}$, which gives the interface shift and the new location of the PCI. However, if the displacement of the PCI over the current time interval is too large, the interface can jump over mesh nodes that were previously not included in the computational subdomain. We have to decide how to treat such new nodes. One approach is described in Section 4.1 of \cite{Gibou2005}, where the authors propose an extrapolation of the numerical temperature solution $T^h$ in normal direction to the interface. The drawback is that this procedure generates a sequence of advection problems \cite{Aslam2004} and therefore adds to the overall complexity of our algorithm. For this reason we refrain from implementing a global extrapolation scheme and we employ an optional adaptive restriction of the time step size. At a given time step $t_n$ we have knowledge about the mesh structure and the advection velocity at any location on the interface, so we define
\begin{equation}
	v_{\text{max}}(t_n) \coloneqq \max \lVert \textbf{v}(\textbf{x}, t_n)\rVert,\hspace{2mm} \textbf{x} \in \Omega.
\end{equation}
Then we choose the next time step such that
\begin{equation}
	t_{n+1} - t_n \leq \frac{h_{\text{min}}}{v_{\text{max}}(t_n)},
\end{equation}
where $h_{\text{min}}$ denotes the minimum element face length.

The second remark concerns the selection of temperature values at the ghost nodes. Going back to Figure \ref{fig:ghost1d}, notice that we assign the phase-change temperature $T_m$ from its exact location on the PCI onto a neighbouring ghost node, which is called \enquote{constant extrapolation} \cite{Gibou2005}. Higher order schemes can be used so that the location where the numerical temperature field satisfies $T^h(\textbf{x},t_n)=T_m$ is shifted closer to the computed position of the PCI. Such schemes are presented in \cite{Gibou2005} within the finite difference framework, but in this work we only focus on the strategy of Section \ref{ghostSplit}.

\subsection{Velocity extension}
\label{velocityExtension}
In Sections \ref{heatFluxReconstruction} and \ref{ghostSplit} we showed how to compute the interface propagation velocity $\textbf{U}(\textbf{x}_{\text{PCI}} , t)$ at the approximate location of the interface intersections. We still need to define the advection velocity in Eq.\ \eqref{eq:levelSet} on all the nodes of the computational mesh. To do so, we use a nearest neighbour classification to decide which value is assigned to each node. That is, we prescribe the Stefan velocity computed at the crossing $C_i$, $i\in\{1,\dots,n_{\text{PCI}}\}$ to all the nodes that are closest to it. Note that we do not require an additional search to find the nearest neighbours of each point, since we already conduct this search in the reinitialization step for the level-set function.

\subsection{Summary of the algorithm}

Now our conceptual workflow, as presented in Fig.\ \ref{fig:2}, can be described in the following algorithmic representation:

\begin{algorithm}[H]
	\SetKwInOut{Input}{Input}
	\SetKwInOut{Output}{Output}
	
	\Input{Initial conditions for velocity, pressure, temperature and the initial geometry of the PCI}
	\Output{Velocity, pressure and temperature fields for both phases and every time step $t_n\in(0,T)$}
	Initialize $\Phi(\textbf{x},0)$ such that $\Phi(\textbf{x}_{\text{PCI}},0)=0$\;
	\Repeat{$t_n<T$}{
	    Update material parameters according to the sign of $\Phi(\textbf{x},t_n)$, see Eq.\ \eqref{eq:properties}\;
	    Solve for the flow and temperature fields at time $t_n$, see governing Eqs.\ \eqref{eq:navierStokes1}, \eqref{eq:navierStokes2}, \eqref{eq:heat}\;
	    Compute the location of the interface $\textbf{x}_{\text{PCI}}$\;
	    Compute the interface velocity $\textbf{U}(\textbf{x}_{\text{PCI}})$ at the crossings with the mesh, see Eq.\ \eqref{eq:meltingVelocity} and Secs.\ \ref{heatFluxReconstruction}, \ref{ghostSplit}\;
	    Extend the interface propagation velocity to all mesh nodes, see Sec.\ \ref{velocityExtension}\;
	    Reinitialize the level-set function\;
	    Solve the level-set Eq.\ \eqref{eq:levelSet}\;
	}
	\Return $\textbf{u}(\textbf{x},t_n)$, $p(\textbf{x},t_n)$, $T(\textbf{x},t_n)$, $\Phi(\textbf{x},t_n)$, $\forall t_n\in(0,T)$\; 
	\caption{Numerical strategy for two-phase melting and solidification problems}
\end{algorithm}

%% file: section4.tex
\section{Numerical results}
\label{section4}

\subsection{Verification: 1D one-phase Stefan problem}

\begin{figure}
    \centering
  \begin{subfigure}[b]{0.55\textwidth}
    \includegraphics[trim=1cm 0.5cm 1cm 0cm,clip,width=0.99\textwidth]{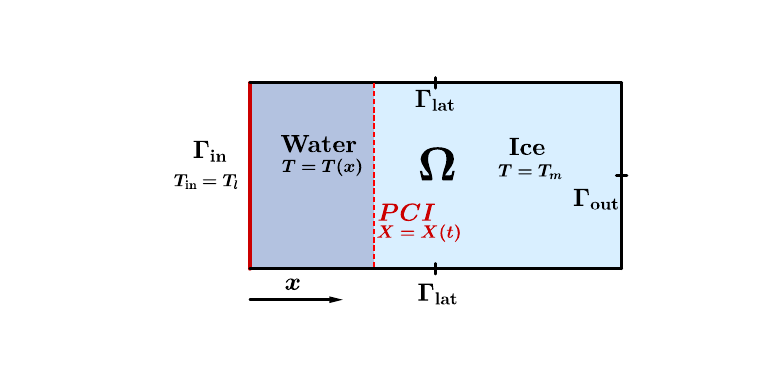}
\subcaption{Schematic of the first test case}
	\label{fig:stefanNumericalSketch}
    \end{subfigure}
    \begin{subfigure}[b]{0.44\textwidth}
    \includegraphics[trim={0cm 0cm 0cm 0cm},clip, width=0.9\textwidth]{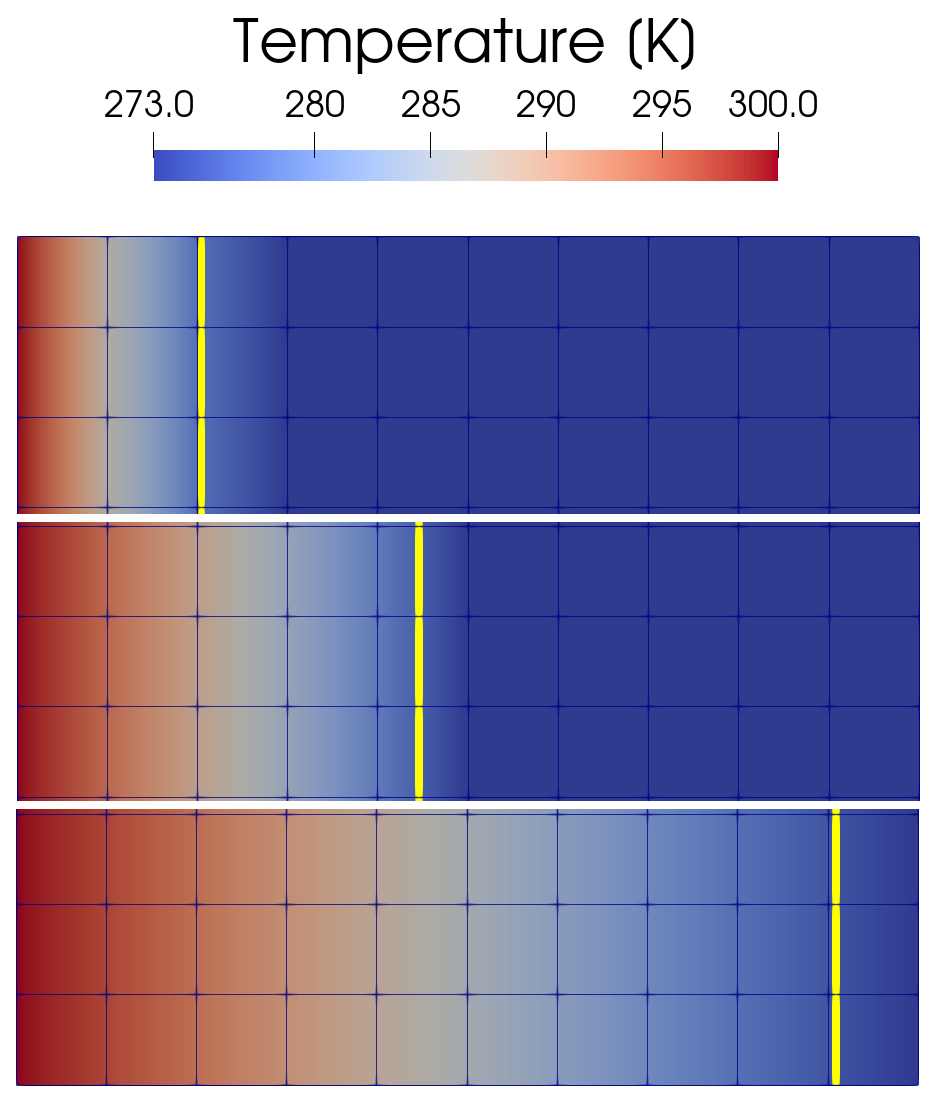}
        \subcaption{Temperature profiles at $t_1=50s$, $t_2=250s$, $t_3=1000s$}
    \label{fig:1dStefanNumericalSolution}
        \end{subfigure}
    \caption{One-phase Stefan problem verification case. On the left, the numerical setup is shown. The ice in initially kept at the constant melting temperature $T_m$. The left boundary $\Gamma_{\text{in}}$ is heated at the temperature $T_{l} > T_m$.
    On the right, the temperature profile at three time instants is shown for a uniform structured grid. For a better visualization we choose a coarse mesh with discretization step $h = $1e-3. The yellow line denotes the PCI.}
 \label{fig:1dStefanTotal}
\end{figure}

\begin{figure}
  \begin{subfigure}[t]{0.45\textwidth}
    \centering
    \includegraphics[trim={0cm 0cm 0cm 0cm},clip, width=0.8\textwidth]{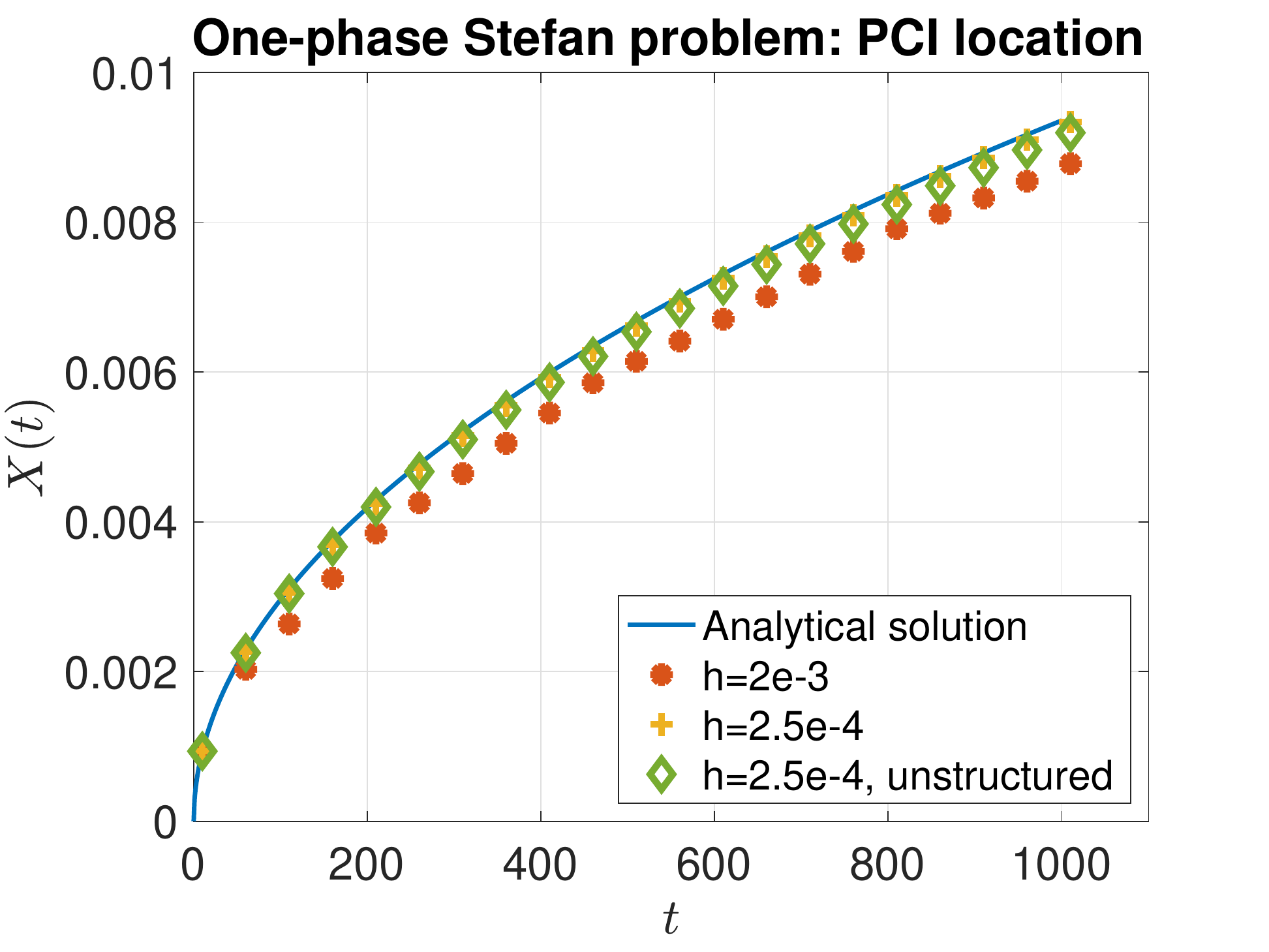}
    \includegraphics[trim={0cm 0cm 0cm 0cm},clip, width=0.8\textwidth]{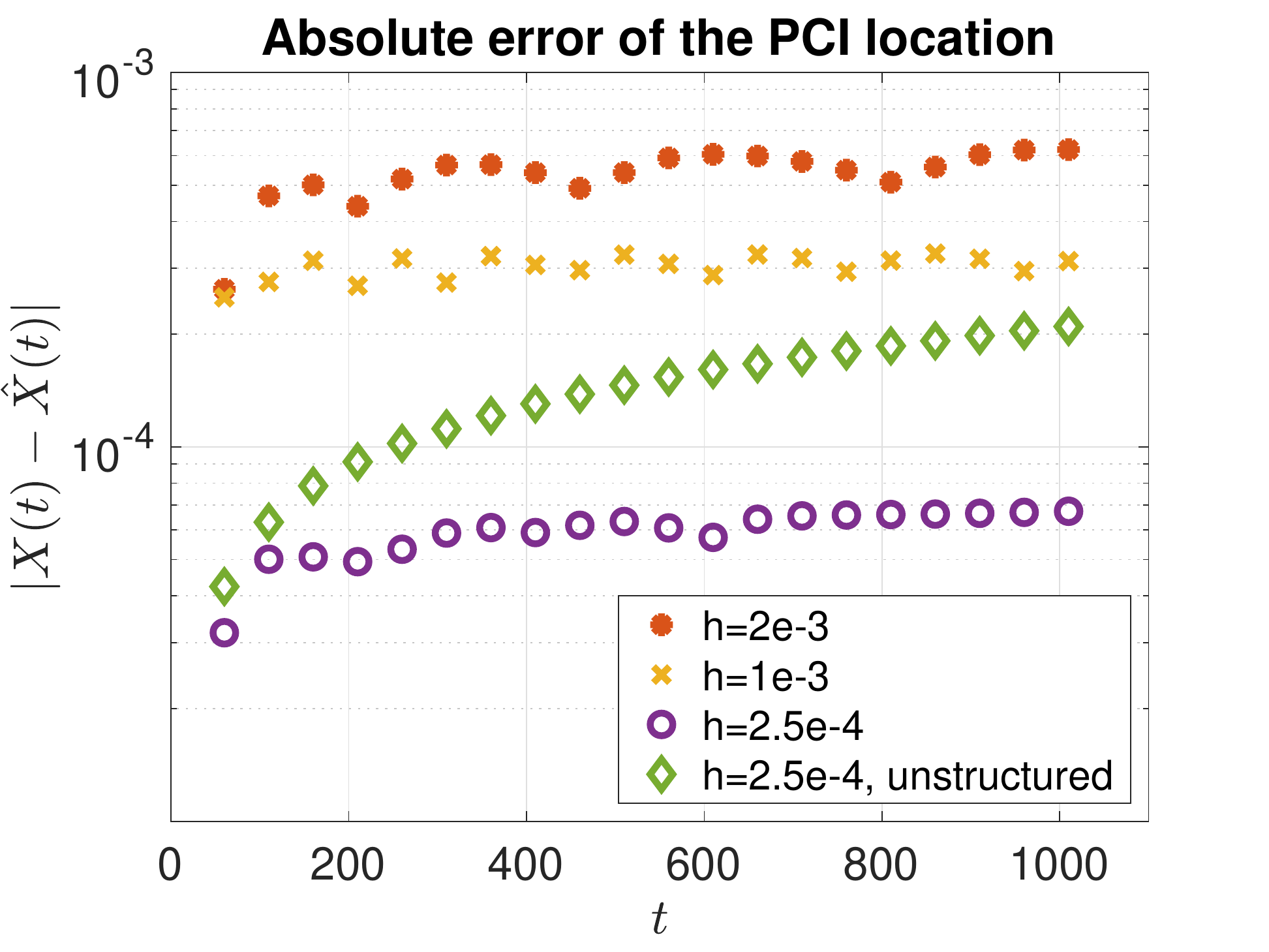}
    \subcaption{Error analysis of the phase-change interface}
             \label{fig:1dStefanInterfaceLocation}
  \end{subfigure}
 \begin{subfigure}[t]{0.45\textwidth}
    \centering
     \includegraphics[trim={0cm 0cm 0cm 0cm},clip, width=0.8\textwidth]{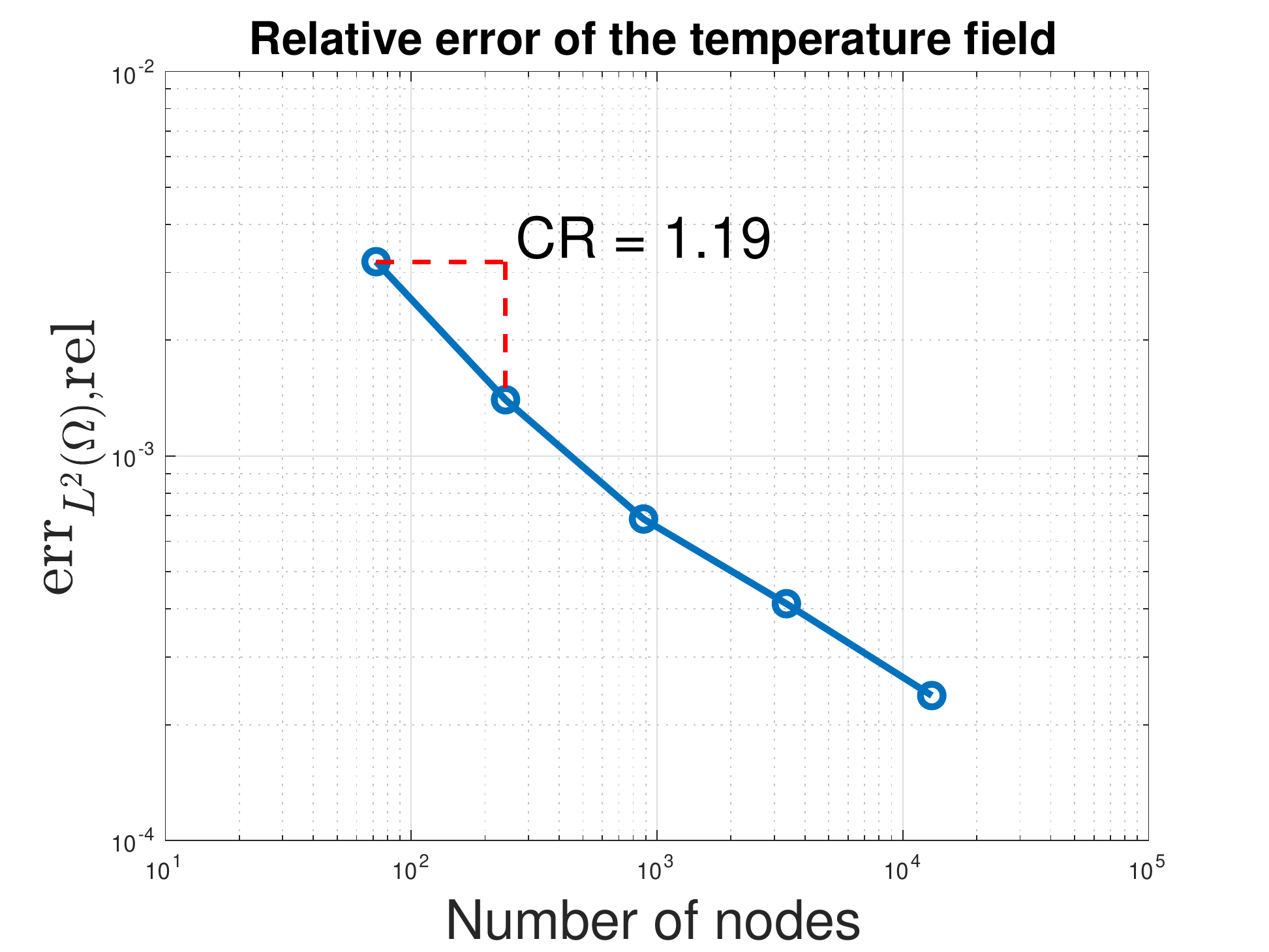}
         \vspace{0.9cm}
   \begin{tabular}{ccc} 
   \toprule
   $h$ & $\text{err}_{L^2(\Omega)}$ & $\text{err}_{L^2(\Omega),\text{rel}}$ \\ 
   \midrule
   2e-3 & 0.9107 & 0.0032 \\
   1e-3 & 0.4116 & 0.0014 \\
   5e-4 & 0.1958 & 6.8551e-4 \\
   2.5e-4 & 0.1179 & 4.1291e-4 \\
   1.25e-4 & 0.0680 & 2.3824e-4 \\
   \bottomrule
	\end{tabular}
	\vspace{0.6cm}
      \subcaption{Error analysis of the temperature profile}
      \label{table:structuredL2error}
    \end{subfigure}
    \caption{One-phase Stefan problem verification case. On the left, the evolution of the PCI is shown for three uniform structured grids of cell size $h$. Additionally, one unstructured grid of triangular elements is considered. The absolute error with respect to the analytical position of the PCI is plotted over time. On the right, the absolute and relative errors of the temperature profile in $L^2$ norm are computed at $t=1000s$ on uniform structured grids. In the top-right figure the relative error is plotted against the number of nodes.}
    \label{fig:stefan1Danalysis}
    
\end{figure}

To verify our method we consider, at first, a 1D one-phase Stefan problem: A slab of ice is initially at constant melting temperature $T_0=T_m$. Then, a constant temperature $T_l>T_m$ is applied at the left boundary $\Gamma_{\text{in}}$, which leads to melting of the ice and causes a phase-change interface (PCI) to propagate to the right, see Fig.\ \ref{fig:stefanNumericalSketch}. Let $\Omega\subset\mathbb{R}^2$ be a bounded domain, let $X(t)$ be the position of the PCI at time $t\in (0,T)$. The governing equations for the problem take the form

\begin{equation}
\begin{aligned}
	\frac{\partial T}{\partial t} = \alpha\frac{\partial^2 T}{\partial x^2} \hspace{5mm} &\textnormal{for}\hspace{2mm}0\leq x < X(t), \hspace{2mm} 0<t<T,\\
    \rho \, h_m \frac{\partial X(t)}{\partial t} = -\kappa \frac{\partial T(X(t),t)}{\partial x} \bigg|_{X^-(t)} \hspace{5mm} &\textnormal{for}\hspace{2mm} 0<t<T, \\
		T = T_l  \hspace{5mm} &\textnormal{on}\hspace{2mm}\Gamma_{\text{in}}, \\
		\nabla T\cdot \textbf{n} = 0 \hspace{5mm} &\textnormal{on}\hspace{2mm}\Gamma_{\text{lat}}\cup\Gamma_{\text{out}}, \\
		T(x,0) = T_m  \hspace{5mm} &\textnormal{in}\hspace{2mm}\Omega,
\end{aligned}
\label{eq:1phaseStefan}
\end{equation}
consisting of the heat equation in the liquid region, in which $\alpha = \kappa / (\rho \, c_p)$ denotes the thermal diffusivity, the Stefan condition and boundary conditions. Note the reduced form of the Stefan condition in comparison to Eq.\ \eqref{eq:meltingVelocity} and the absence of a heat equation for the solid phase. Both are due to the constant temperature (hence zero temperature gradient) in the solid phase, which gives rise to the notion of the one-phase Stefan problem. Thus, the terminology \enquote{one-phase} acknowledges the fact that we solve for the temperature in the liquid portion only.
The analytical solution of Problem \eqref{eq:1phaseStefan} can be found using a similarity approach. Following \cite{Jonsson2013}, we obtain
\begin{equation}
\begin{aligned}
	\hat{T}(x,t) = T_l - (T_l - T_m)\frac{\erf\left(\frac{x}{2\sqrt{\alpha t}}\right)}{\erf\left(\frac{X(t)}{2\sqrt{\alpha t}}\right)} \hspace{5mm} &\textnormal{for}\hspace{2mm}0\leq x < X(t), \hspace{2mm} 0<t<T,\\
   \hat{X}(t) = 2\sqrt{\alpha} \, \lambda\sqrt{t} \hspace{5mm} &\textnormal{for}\hspace{2mm} 0<t<T, \\
\end{aligned}
\label{eq:1DstefanAnalytical}
\end{equation}
where $\lambda$ is the unique root of the monotonic function 
\begin{equation}
	f(\lambda)\coloneqq \frac{c_p(T_m - T_l)}{h_m}\frac{e^{-\lambda^2}}{\sqrt{\pi}\erf(\lambda)} - \lambda,
\end{equation}
and $\erf$ denotes the error function $\erf(x) \coloneqq 2/\pi \int_0^x e^{-y^2}\diff y$. Note that the term \enquote{1D} highlights the dependence of both the temperature and the PCI location on the sole $x$ coordinate, but the problem has a bidimensional setting. For the numerical simulation we consider a $0.01\times0.01$ square domain and compute 2000 time steps with $\Delta t = 0.5s$. The initial time is $t_0=10s$, where we prescribe the analytical solutions for the temperature and the PCI location, see Eq.\ \eqref{eq:1DstefanAnalytical}, as initial conditions. The physical parameters are selected according to water ice, namely $T_l = \SI{300}{\kelvin}$, $T_m = \SI{273}{\kelvin}$, $h_m = \SI{333700}{\joule\per\kg}$, $\rho = \SI{1000}{\kg\per\m\cubed}$, $c_p = \SI{4200}{\joule\per\kg\per\kelvin}$, $\kappa = \SI{0.6}{\watt\per\m\per\kelvin}$.\\
Figure \ref{fig:1dStefanNumericalSolution} shows the evolving temperature profile for one particular simulation setting, that is a uniform structured grid with spatial cell size $h = $1e-3. It is clearly visible that the phase interface propagates from left to right for increasing simulation times. At each time step, the temperature gradually decreases from the left Dirichlet boundary towards the melting temperature $T_m$ at the interface, while it stays constant in the solid phase.\\
Figure \ref{fig:1dStefanInterfaceLocation} shows the evolving location of the PCI for a structured grid of different cell sizes $h$ as well as for an unstructured grid with triangular elements. For visualization purposes, we plot the numerical values every 100 time steps. The PCI positional error is computed against the analytically predicted PCI location from Eq.\ \eqref{eq:1DstefanAnalytical}. It can be seen that such error slightly increases with time and diminishes as the grid gets finer. Likewise, the error introduced by the ghost split decreases for $h \rightarrow 0$, as the melting temperature is imposed on a node closer to the exact position of the PCI. Finally, Table \ref{table:structuredL2error} displays the $L^2$ error on a cascade of refined structured meshes with respect to the analytical value of the temperature field at $t=1000s$, see Eq.\ \eqref{eq:1DstefanAnalytical}. The same values are plotted against the number of nodes, which results in a convergence rate (CR) of 1.19. This is lower than the expected second order associated with
the employed FEM discretization of the heat equation alone. Recall, however, that here the heat equation is coupled to the interface propagation, and that assigning the melting temperature at the ghost nodes introduces an additional
error. Higher order schemes for the temperature extrapolation are available, see Section \ref{temperatureExtrapolation}, and can be investigated
in the future.

\subsection{Phase-change coupled 2D lid-driven cavity problem}
\begin{figure}
\begin{subfigure}[b]{0.49\textwidth}
	\centering
	\includegraphics[trim=0cm 0cm 0cm 0cm,clip,width=0.7\textwidth]{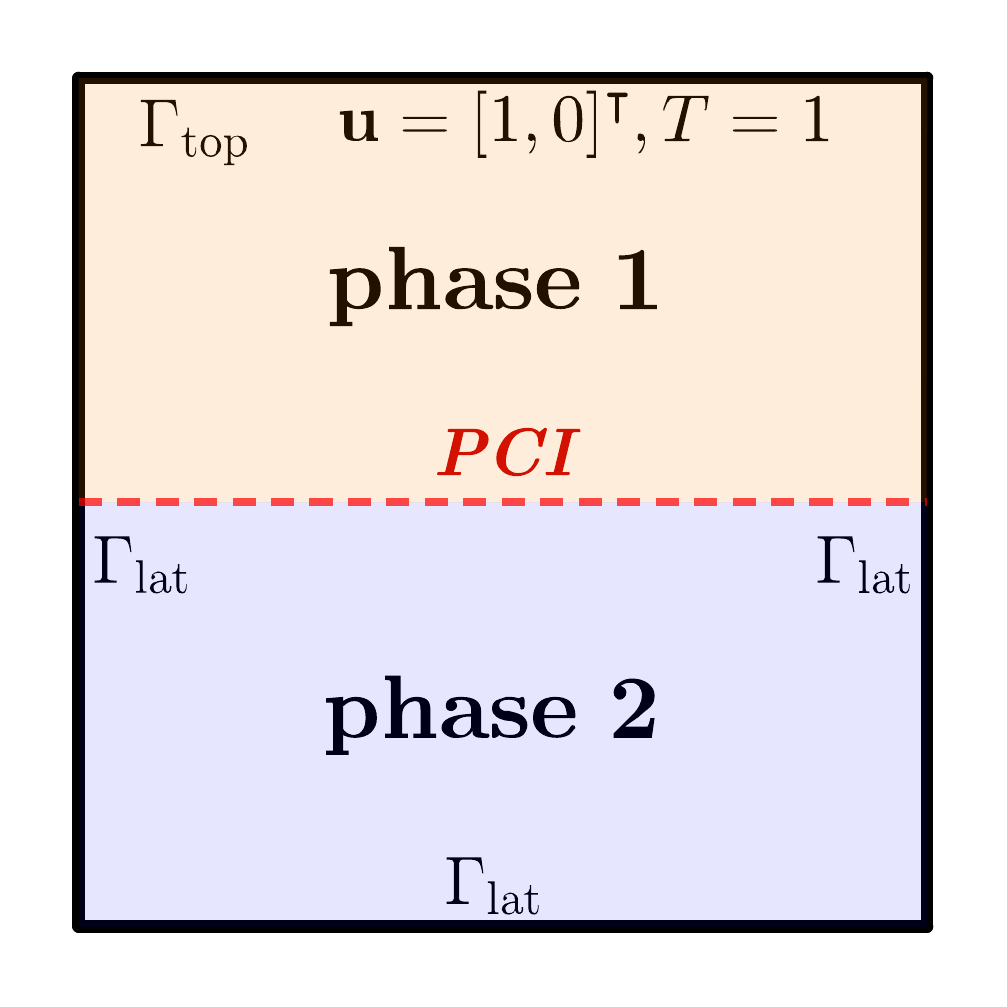}
	\subcaption{Schematic of the second test case}
	\label{fig:lidDrivenDomain}
\end{subfigure}
	\begin{subfigure}[b]{0.49\textwidth}
\begin{center}
\begin{tabular}{lcc} 
   \toprule
   \textbf{Parameter} & \textbf{phase 1} & \textbf{phase 2} \\ 
   \midrule
   $\rho$ $[\si{\kg\per\m\cubed}]$ & 2 & 1 \\
   $c_p$ $[\si{\joule\per\kg\per\kelvin}]$ & 1e3 & 1 \\
   $\kappa$ $[\si{\watt\per\m\per\kelvin}]$ & 1 & 1 \\
   $\mu$ $[\si{\kg\per\m\per\s}]$ & 1 & 1e4 \\
   \midrule
   $h_m$ $[\si{\joule\per\kg}]$ & 1 & - \\
   \bottomrule
\end{tabular}
\end{center}
\subcaption{Material properties}
\label{table:lidDrivenMaterials}
\end{subfigure}
\caption{Setup of the phase-change coupled 2D lid-driven cavity problem. On the left, the computational domain is presented. Initially, the top half is covered by the liquid phase, while the bottom half is covered by the solid phase. On lateral and bottom boundaries $\Gamma_{\text{lat}}$ we impose homogeneous Dirichlet boundary conditions for the velocity field and homogeneous Neumann conditions for the temperature field. Constant temperature and velocity in the $x$ direction are prescribed on the top edge $\Gamma_{\text{top}}$. On the right, we list the material properties.}
\label{fig:lidDrivenTotal}
\end{figure}

\begin{figure}

    \centering
  \begin{subfigure}[c]{0.3\textwidth}
    \includegraphics[trim={0.5cm 0.7cm 0.5cm 0.7cm},clip, width=\textwidth]{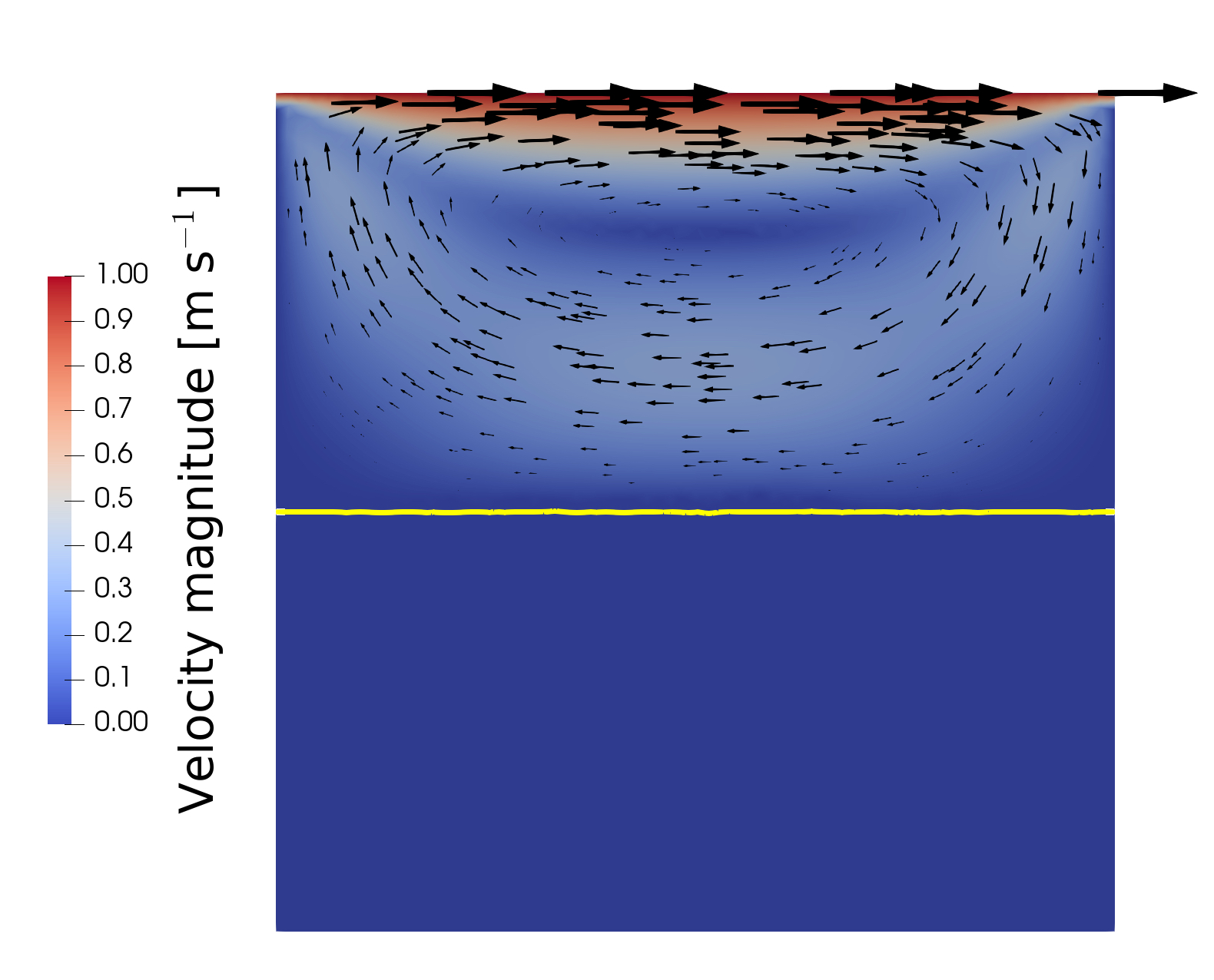}
     \includegraphics[trim={0.5cm 0.7cm 0.5cm 0.7cm},clip, width=\textwidth]{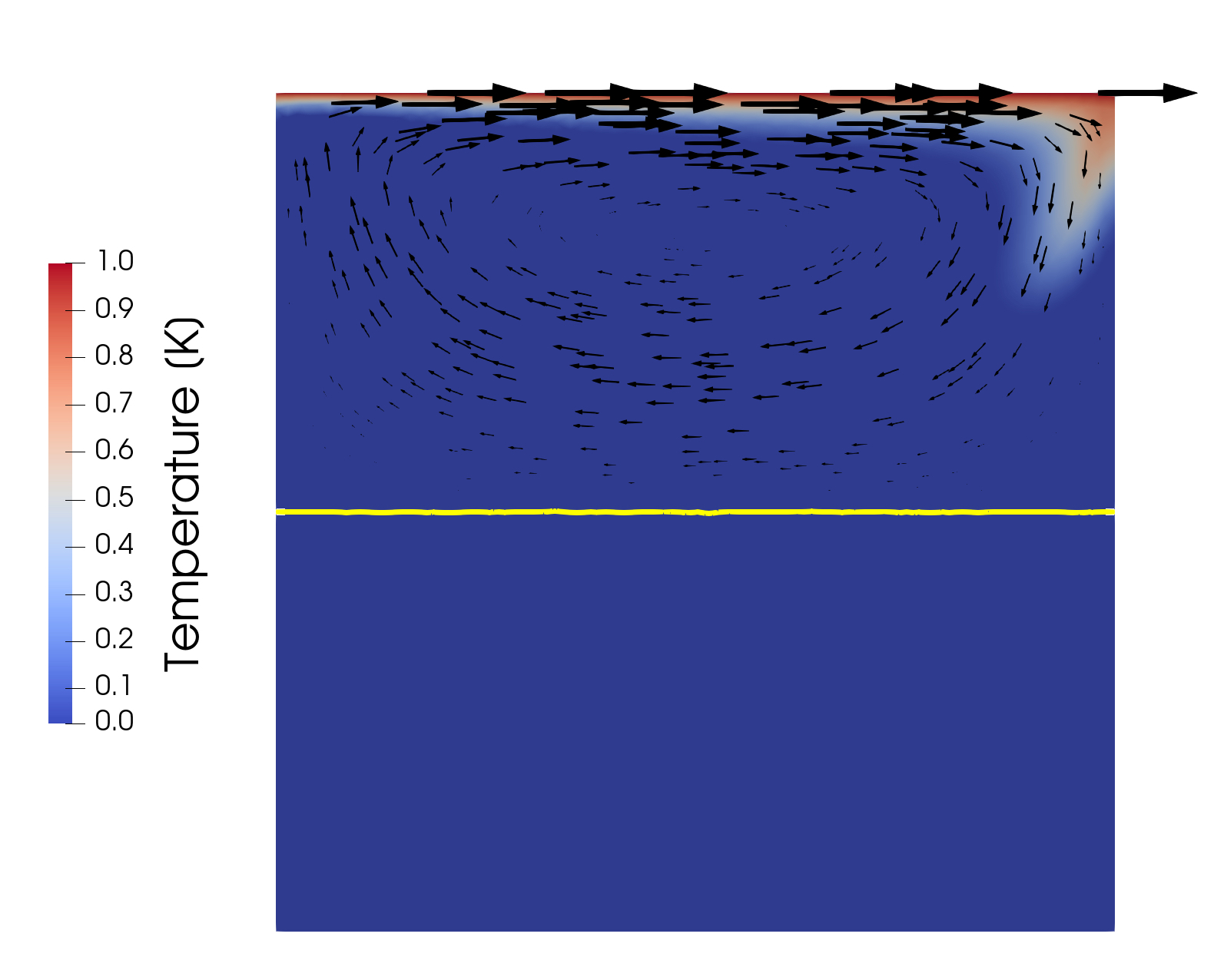}
        \subcaption{$t=1s$}
         \label{fig:lidDrivenVelocityT10}
    \end{subfigure}
    \begin{subfigure}[c]{0.3\textwidth}
    \includegraphics[trim={0.5cm 0.7cm 0.5cm 0.7cm},clip, width=\textwidth]{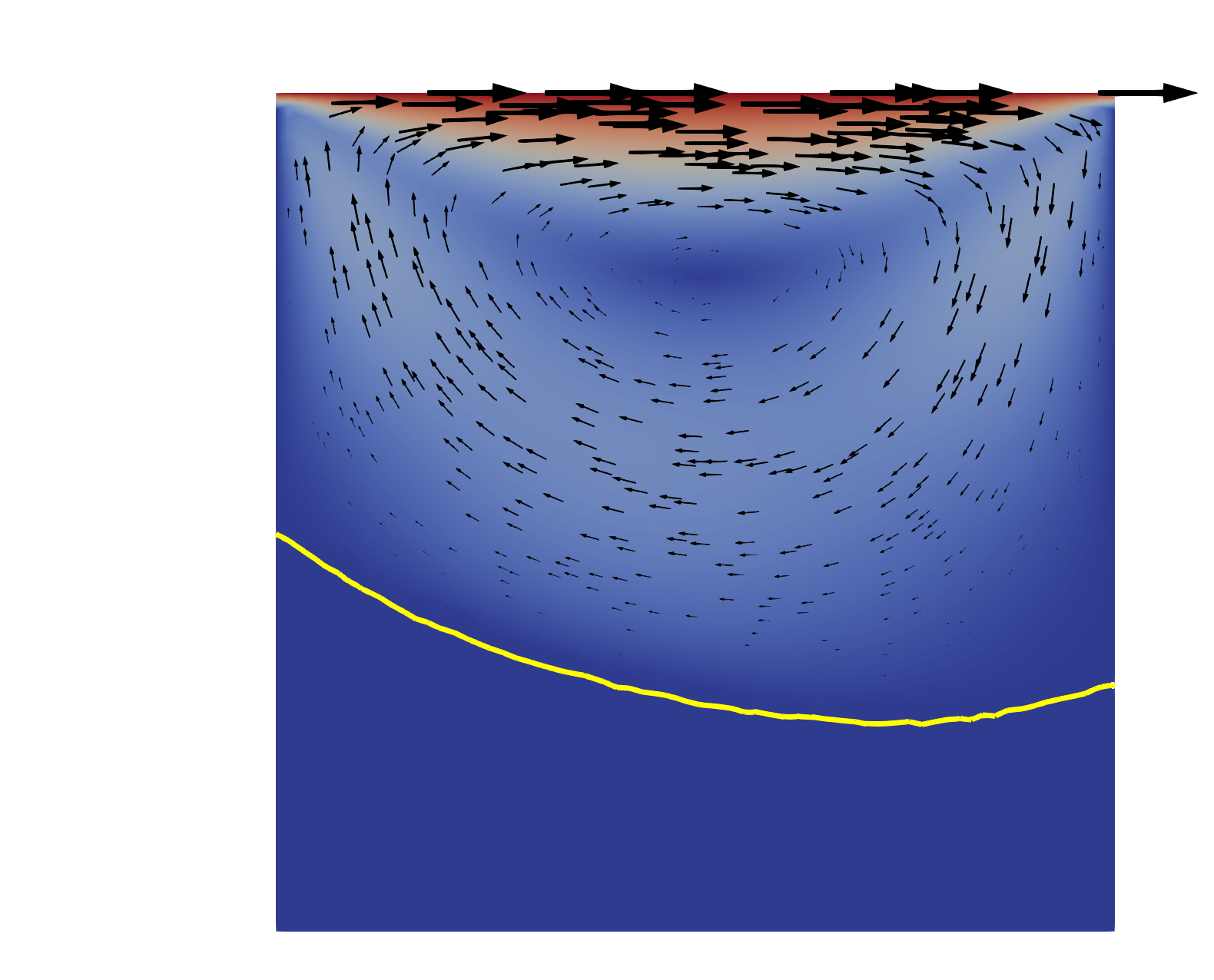}
       \includegraphics[trim={0.5cm 0.7cm 0.5cm 0.7cm},clip, width=\textwidth]{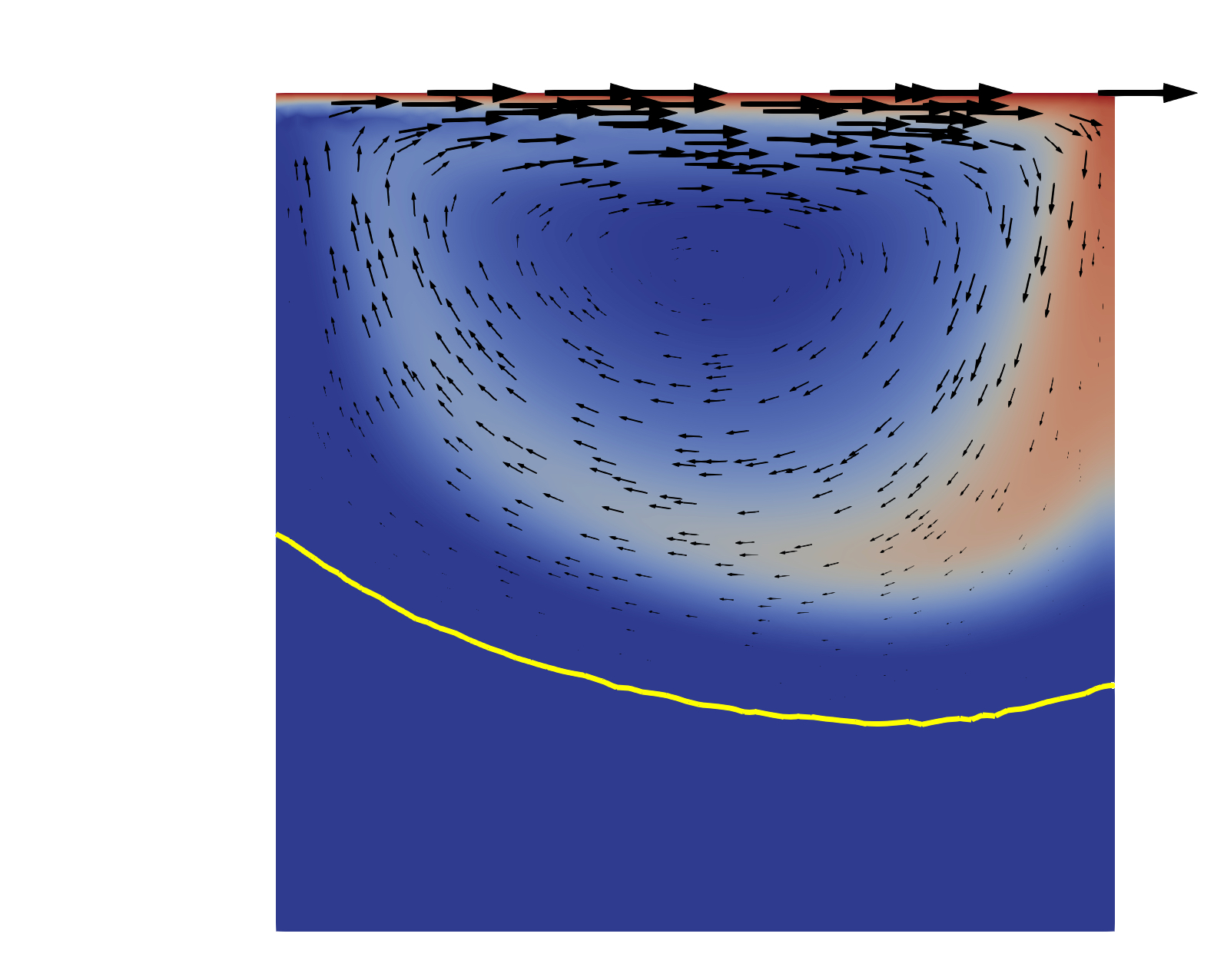}
        \subcaption{$t=10s$}
    \end{subfigure}
    \begin{subfigure}[c]{0.3\textwidth}
    \includegraphics[trim={0.5cm 0.7cm 0.5cm 0.7cm},clip, width=\textwidth]{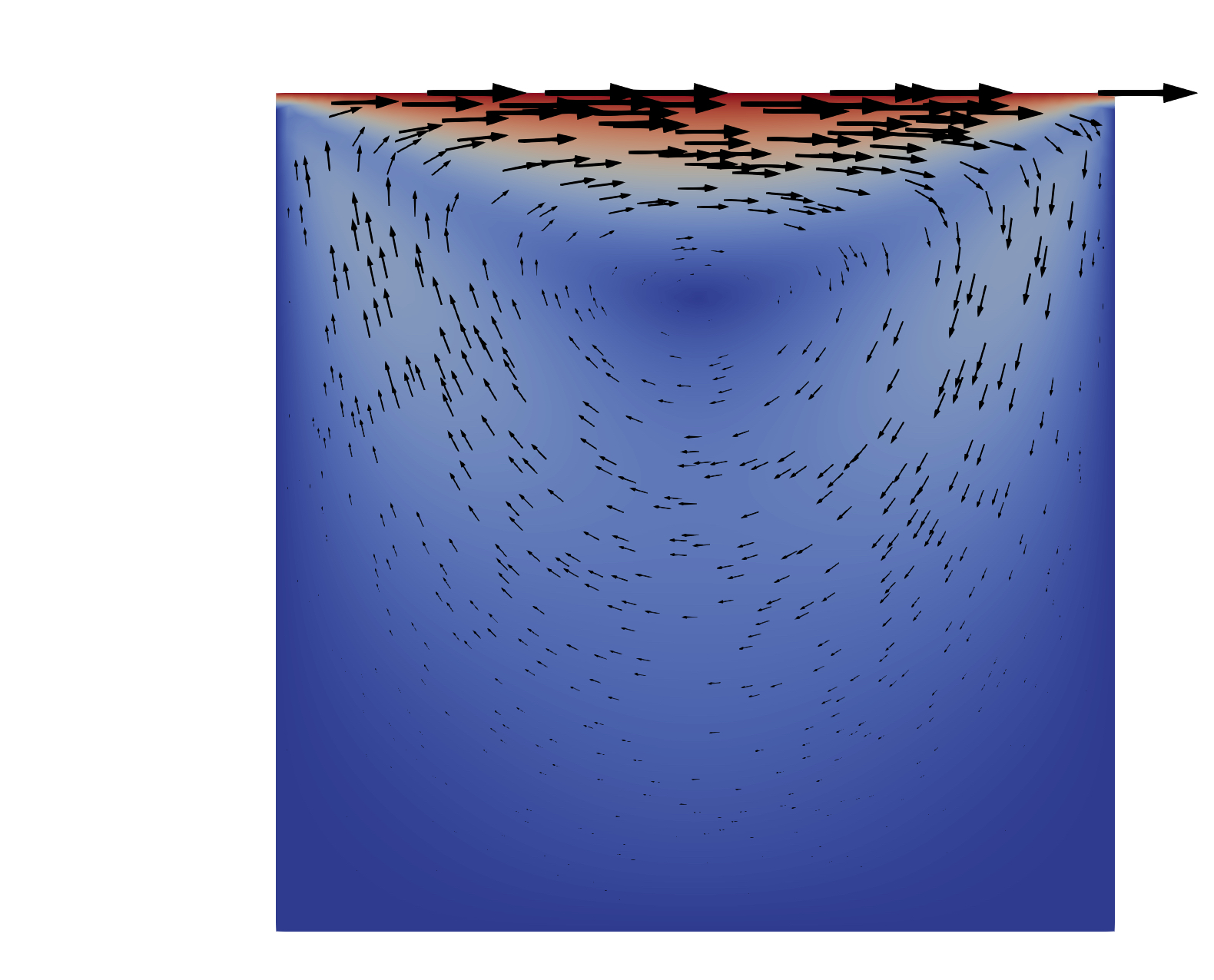}
      \includegraphics[trim={0.5cm 0.7cm 0.5cm 0.7cm},clip, width=\textwidth]{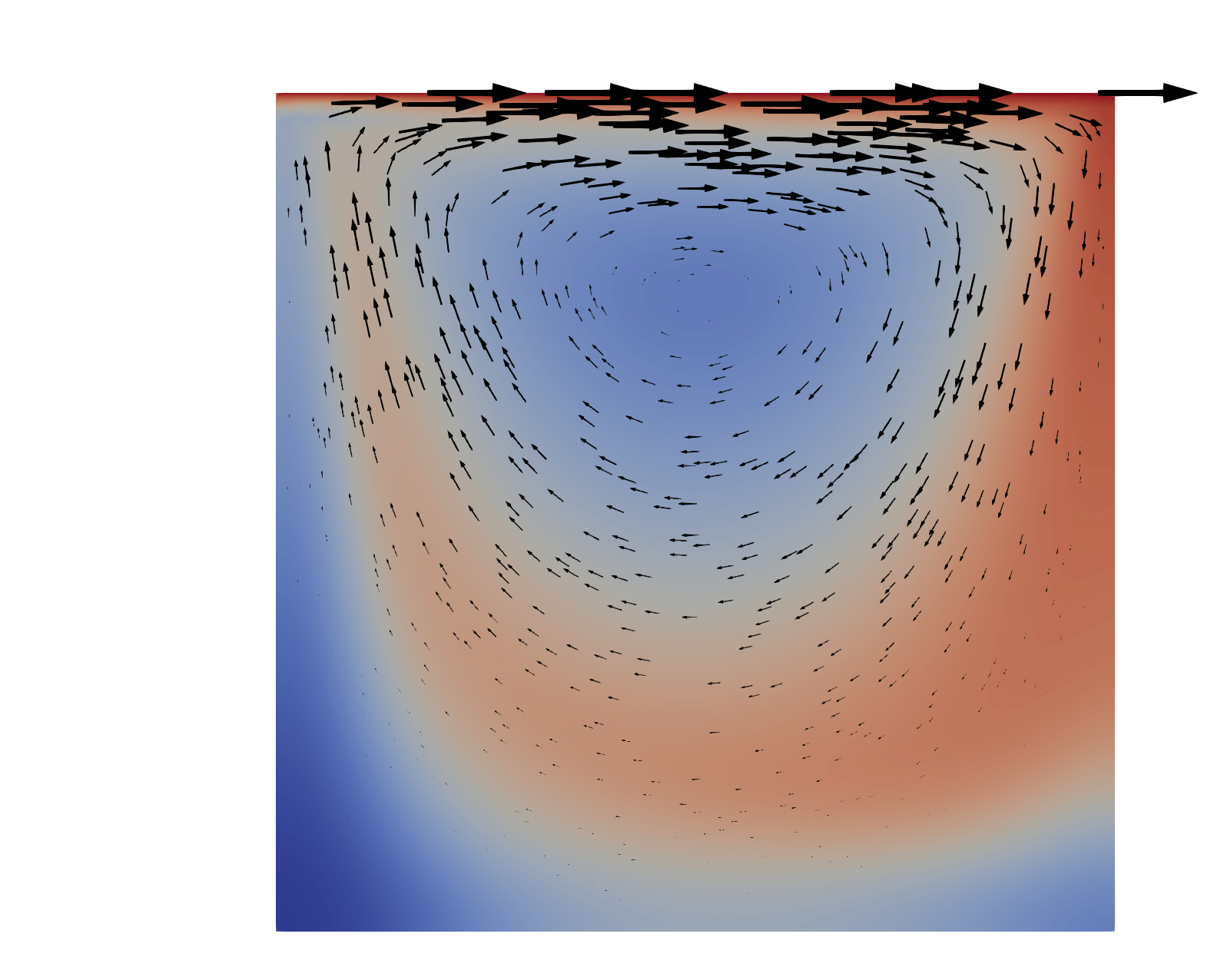}
        \subcaption{$t=40s$}
        \label{fig:lidDrivenVelocityT400}
    \end{subfigure}
    
    \caption{Phase-change coupled 2D lid-driven cavity problem. The velocity magnitude (top) and the temperature field (bottom) are shown at three time instants. Black arrows represent the velocity vectors at each point, their size is proportional to the velocity magnitude. The yellow line denotes the phase-change interface.}
    \label{fig:lidDrivenCavityVelocity}
    
\end{figure}

For the second test case we move to a problem that shows the complete workflow presented in Figure \ref{fig:2}. In particular, we solve both for the evolving velocity and temperature fields and examine different material properties for the two phases. Let us consider the $1\times1$ domain in Fig.\ \ref{fig:lidDrivenDomain}, in which the top half, i.e.\ for $y>0.5$, initially is in liquid state, while the bottom half is in solid state. At the lateral and bottom boundaries, indicated by $\Gamma_{\text{lat}}$, we prescribe homogeneous Dirichlet boundary conditions for the velocity and homogeneous Neumann conditions for temperature. At the top edge $\Gamma_{\text{top}}$ we impose a constant temperature $T=1$ and constant velocity in $x$ direction $\textbf{u}=[1,0]^\intercal$. The initial temperature is $T_m=0$ over the whole domain. The material properties are shown in Table \ref{table:lidDrivenMaterials}. They are purely fictional, yet we choose them in order to emphasize the role of convection in the heat equation. Note that the solid is not modeled explicitly, but instead as a fluid with a relatively high viscosity.
We solve the problem on an unstructured grid with cell size $h=0.02$. We simulate 500 time steps with $\Delta t = 0.1s$. Figure \ref{fig:lidDrivenCavityVelocity} shows the computed flow field at three different time instants. Right after the start, we retrieve the expected clockwise circulation of the fluid in the top half, i.e.\ the liquid domain, as shown by the black velocity vectors, see Fig.\ \ref{fig:lidDrivenVelocityT10}. After 400 time steps, the solid material is completely molten and the liquid can circulate in the whole domain. This gives the familiar lid driven cavity flow in a square domain \cite{kuhlmann2018}. Figure \ref{fig:lidDrivenCavityVelocity} shows the temperature profiles at the same time instants. Note that the temperature evolution is driven by the convection of the flow field. For that reason we can observe that the right side of the domain melts faster.

\begin{figure}

    \centering
  \begin{subfigure}[c]{0.49\textwidth}
    \includegraphics[trim={0.5cm 0.7cm 0.5cm 0.7cm},clip, width=0.9\textwidth]{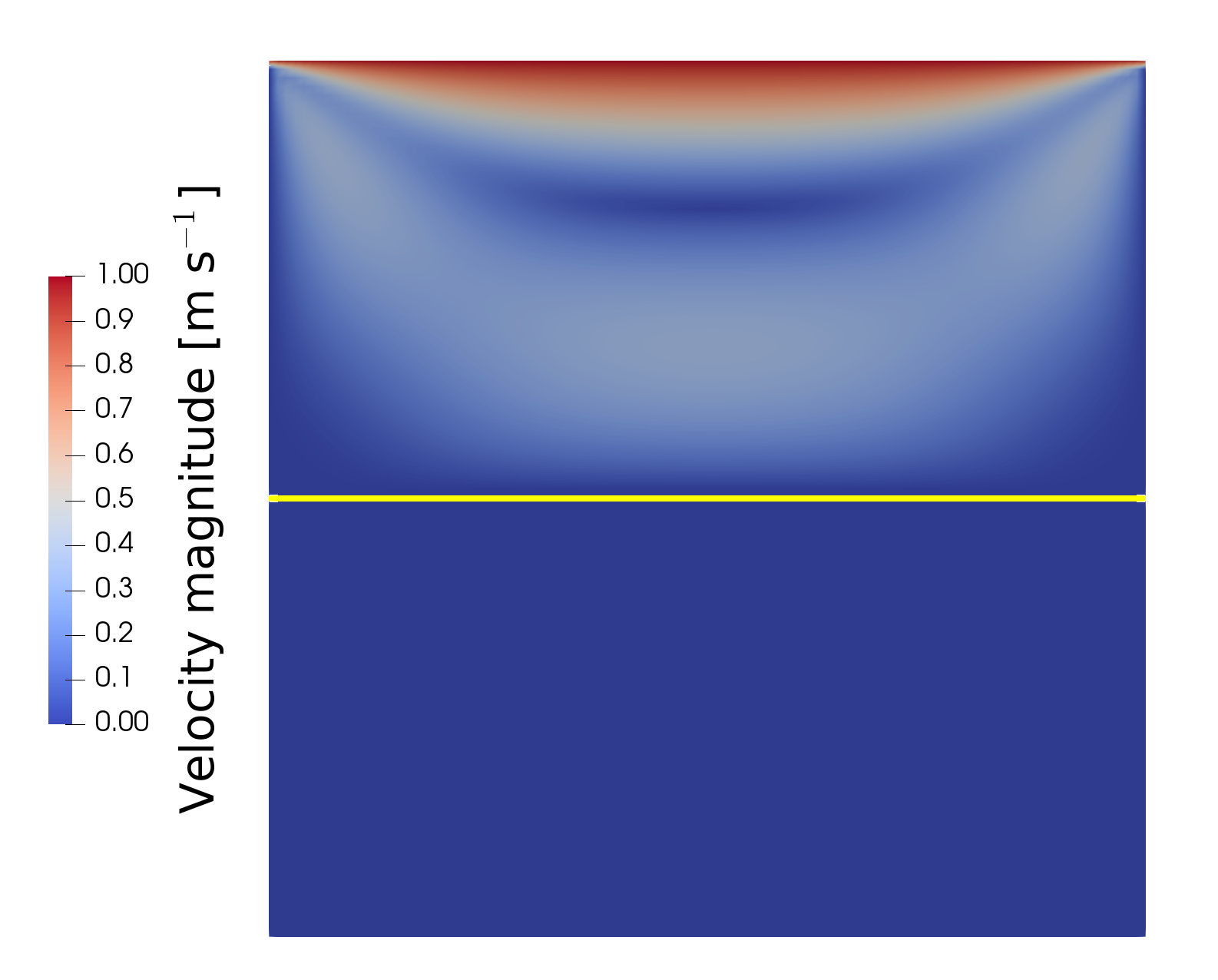}
    \caption{Full simulation for coupled flow and temperature}
    \end{subfigure}
    \begin{subfigure}[c]{0.49\textwidth}
    \includegraphics[trim={0.5cm 0.7cm 0.5cm 0.7cm},clip, width=0.9\textwidth]{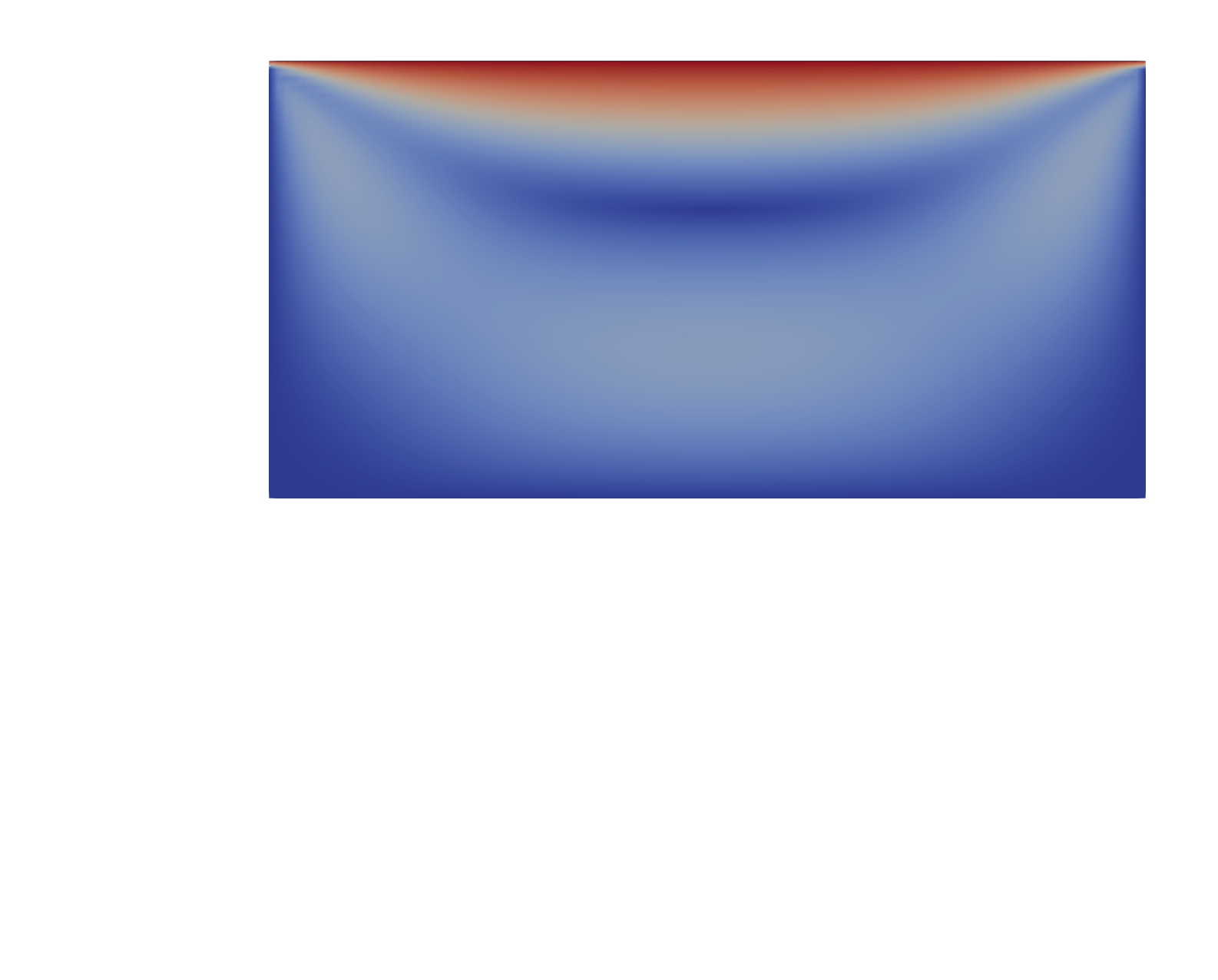}
    \caption{Flow field simulation on the top half domain}
    \end{subfigure}
    
    \caption{Phase-change coupled 2D lid-driven cavity problem. We show the velocity magnitude at $t=1s$. Both meshes have a uniform discretization step $h=0.01$. On the left, the result from the coupled phase-change simulation of the melting problem, where the yellow line denotes the PCI. On the right, the result from the classical velocity simulation on the liquid region. Note that there are no visual differences between the two solutions.}
    \label{fig:lidDrivenComparisonU}
    
\end{figure}

Notice that the simulation we have just described embeds a conventional lid-driven cavity problem with no-slip, hence homogeneous Dirichlet conditions, at the boundary defined by the PCI. In order to demonstrate this explicitly, we compare both setups in a numerical experiment. First, we simulate the phase-change coupled lid-driven cavity scenario from before. For this scenario, however, we compute only 10 time steps, so that the PCI has not yet changed its position at the end of the simulation. Second, we compute the velocity profile for a classical lid-driven cavity problem, where we impose a no-slip boundary condition at the bottom edge. This time the computational domain corresponds to the liquid region of the first simulation. Both problems are computed on a uniform structured grid of cell size $h=0.01$. 
Figure \ref{fig:lidDrivenComparisonU} shows the velocity magnitude at the final time instant $t=1s$. Recall, that in the phase-change coupled simulation there is no zero-velocity imposition at the PCI, but the velocity relaxation towards the interface is rather a result of selecting an extremely high value of the viscosity in the solid phase. Despite that, there are no visual differences with the second plot, where we have set a no-slip condition at the bottom. The relative errors between both simulation runs are $\text{err}^u_{L^2(\Omega),\text{rel}} = 0.0282$ and $\text{err}^v_{L^2(\Omega),\text{rel}} = 0.0322$ for the two components $u$ and $v$. We recall that our implementation considers a PCI of fixed thickness, so that the material properties can vary smoothly for the sake of numerical stability, see Eq.\ \ref{eq:heaviside}. The errors above have been computed with the value $\epsilon=0.001$. If we select a sharp interface of $\epsilon=0$ instead, the errors drop to $\text{err}^u_{L^2(\Omega),\text{rel}} = 2.6037$e-4 and $\text{err}^v_{L^2(\Omega),\text{rel}} = 6.3687$e-4. This shows that we are able to represent the zero-velocity boundary condition at an immersed phase-change boundary condition.

\subsection{Phase-change coupled 2D corner flow - interface and temperature evolution}
\label{2DmeltingProbe}
The last example is inspired by one of our applications of interest, namely the spatio-temporal evolution of a melt channel that develops as a thermal melting cryorobot descents into the ice \cite{dachwald2014}. While we gained significant knowledge on the melting probe performance in recent years \cite{schuller2016curvilinear,schuller2019melting}, a holistic model that integrates the cryobot dynamics with thermo-fluidmechanically coupled processes in the melt channel is missing to date. In particular, we need the capability to model convection-coupled phase-change in complex corner flow geometries.

\begin{figure}
\centering
\includegraphics[trim=0cm 1.5cm 0cm 1cm,clip,width=0.95\textwidth]{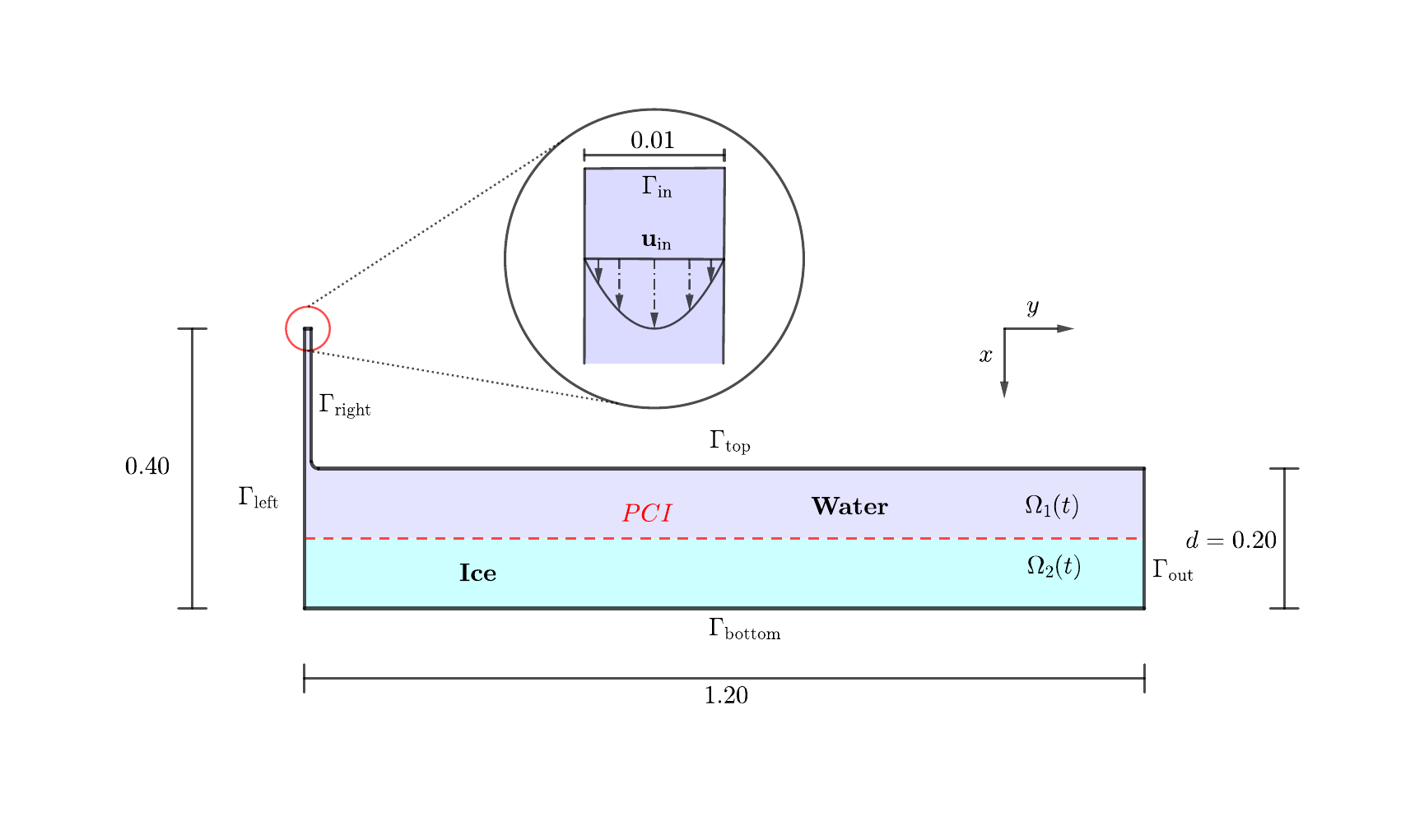}
\caption{Phase-change coupled 2D corner flow. This setup resembles one half of an idealized melting probe moving to the left. The liquid region $\Omega_2(t)$ is shown in blue. We consider a parabolic velocity profile at the inflow $\Gamma_{\text{in}}$ and a constant temperature on $\Gamma_{\text{right}}$, which causes the ice block to melt.}
\label{fig:2dProbeDomain}
\end{figure}

\begin{figure}

    \centering
  \begin{subfigure}[c]{0.7\textwidth}
    \includegraphics[trim={0cm 0cm 0cm 0cm},clip, width=\textwidth]{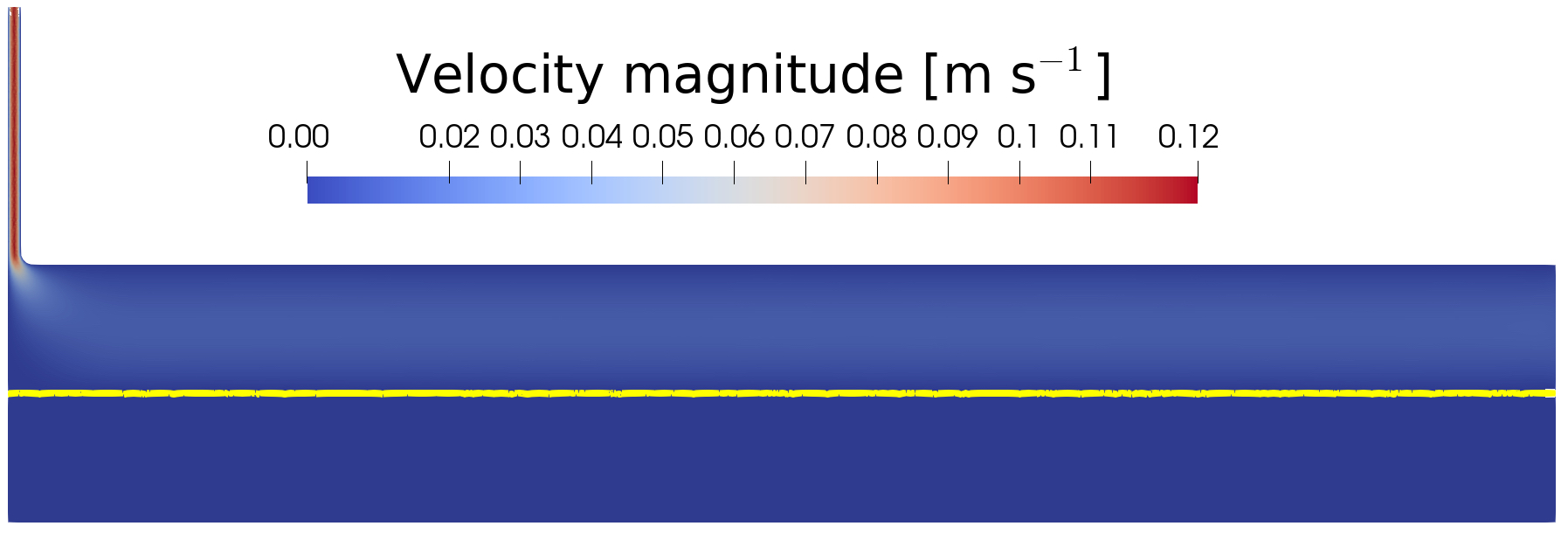}
    \caption{$t = 10s$}
    \end{subfigure}
    \begin{subfigure}[c]{0.7\textwidth}
    \includegraphics[trim={0cm 0cm 0cm 0cm},clip, width=\textwidth]{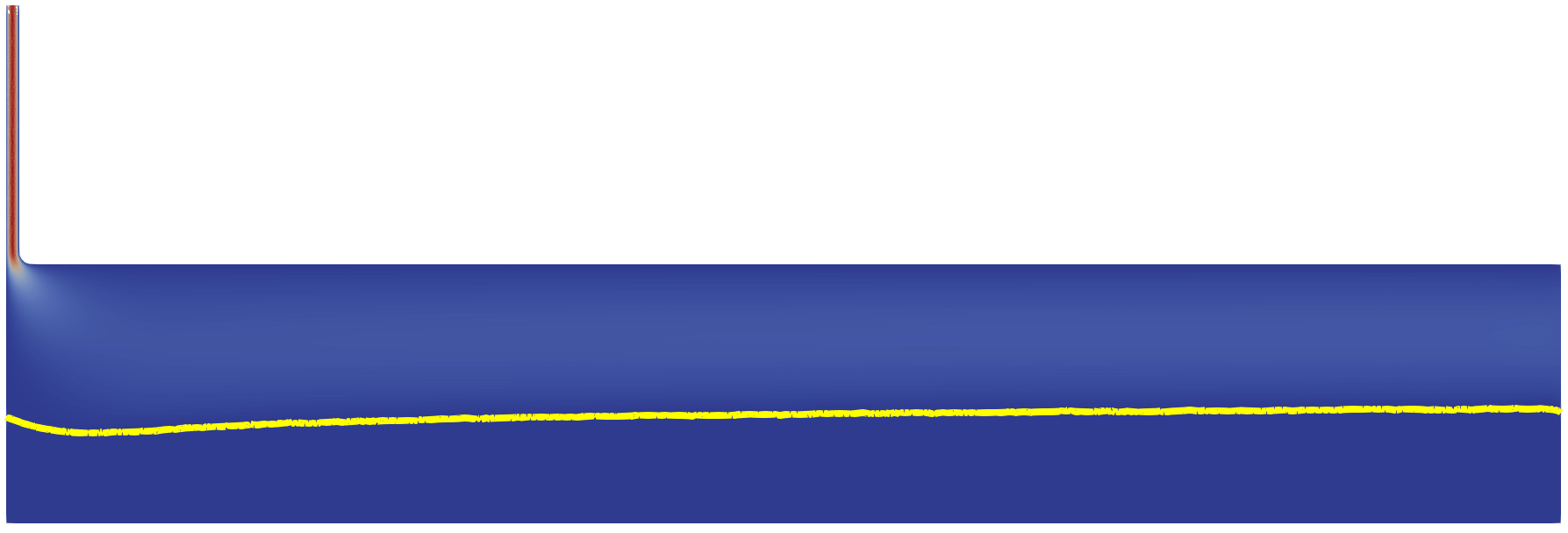}
    \caption{$t = 2500s$}
    \end{subfigure}
    
    \caption{Phase-change coupled 2D corner flow. The velocity magnitude is shown at two different time instants. The yellow line denotes the phase-change interface. As expected, the left side of the solid region melts faster due to the heat transported by the water. }
    \label{fig:probe2Dvelocity}
    
\end{figure}

\begin{figure}

    \centering
  \begin{subfigure}[b]{0.49\textwidth}
    \includegraphics[trim={0cm 0cm 0cm 0cm},clip, width=\textwidth]{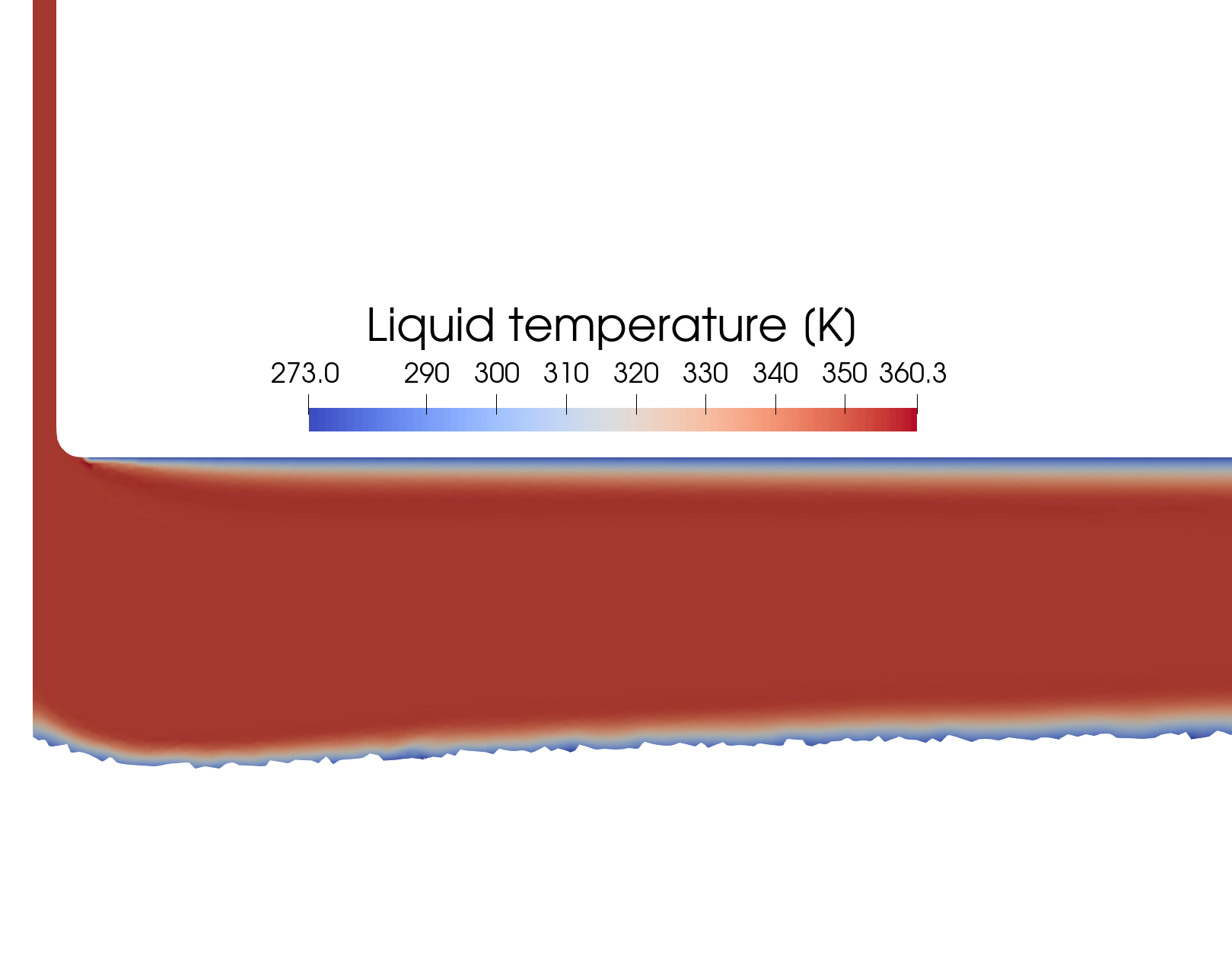}
    \caption{liquid phase}
    \end{subfigure}
    \begin{subfigure}[b]{0.49\textwidth}
    \includegraphics[trim={0cm 0cm 0cm 0cm},clip, width=\textwidth]{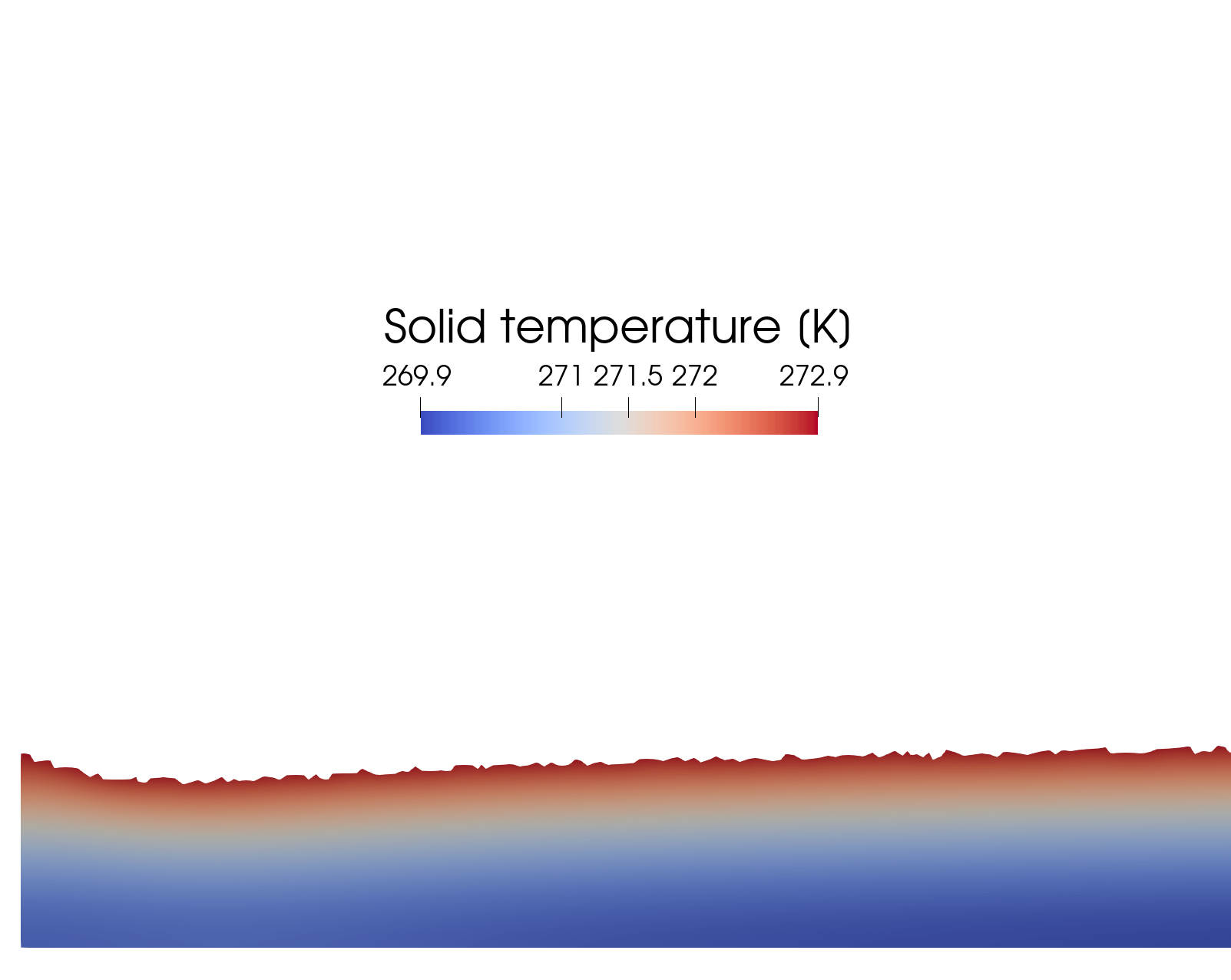}
    \caption{solid phase}
    \end{subfigure}
    
    \caption{Phase-change coupled 2D corner flow. The two temperature profiles for the liquid and solid phases are shown at $t=2500s$. For a better visualization we display only the left part of the domain, that is for $y<0.6$. We recall that due to the ghost-cell method we solve two separate temperature problems for each phase. The two subdomains are divided by the phase-change interface.}
    \label{fig:probe2DTemperature}
    
\end{figure}

Figure \ref{fig:2dProbeDomain} represents the geometry of our final test case: The inflow into a small channel is diverted by 90 degrees into a wider outflow channel. The inflow channel has fixed boundaries, while the wider outflow channel contains two different phases separated by an evolving PCI. Note, that this setup resembles one half of an idealized melting probe moving to the left. The inflow channel has a width of $1/20$ with respect to the outflow channel. 
We impose a parabolic velocity profile at the inflow $\Gamma_{\text{in}}$. Furthermore, we impose zero-velocity conditions at each boundary except for the inflow and the outflow boundaries. We have Dirichlet temperature conditions at the interior boundaries and at the inflow boundary, $T=353$ and $T=278$ respectively, and homogeneous Neumann conditions everywhere else. A summary of all boundary and initial conditions for the problem is
\begin{equation}
	\begin{aligned}
		\textbf{u}(\textbf{x},0) = \textbf{0} \hspace{5mm} &\textnormal{in}\hspace{2mm}\Omega_{1,0} \cup \Omega_{2,0},\\
		T(\textbf{x},0) = 273 \hspace{5mm} &\textnormal{in}\hspace{2mm}\Omega_{1,0},\\
		T(\textbf{x},0) = 268 \hspace{5mm} &\textnormal{in}\hspace{2mm}\Omega_{2,0},\\
		\textbf{u}=(5000y \, (0.01-y),0)^\intercal  \hspace{5mm} &\textnormal{on}\hspace{2mm}\Gamma_{\text{in}}, \\
		\textbf{u}=\textbf{0}  \hspace{5mm} &\textnormal{on}\hspace{2mm} \Gamma_{\text{left}} \cup \Gamma_{\text{right}} \cup \Gamma_{\text{top}} \cup \Gamma_{\text{bottom}}, \\
		T = 353  \hspace{5mm} &\textnormal{on}\hspace{2mm}\Gamma_{\text{in}} \cup \Gamma_{\text{right}}, \\
		T = 278  \hspace{5mm} &\textnormal{on}\hspace{2mm}\Gamma_{\text{top}}.
	\end{aligned}
\label{eq:probe2dConditions}
\end{equation}
All values are in SI units. The initial location of the PCI is at $x=0.3$. The material properties associated to the two phases are shown in Table \ref{table:2dProbeMaterials}. Note that different values are assigned to each material parameter and to each phase in the system, so that we can more realistically replicate the behaviour of water and ice. We simulate the problem on an unstructured grid with 49196 nodes and compute 500 time steps with $\Delta t=5s$. Figure \ref{fig:probe2Dvelocity} shows the velocity profile at two different time instances. It is clearly visible how a bulge forms in the wider part of the channel due to the warm diverted inflow.
Figure \ref{fig:probe2DTemperature} displays the temperature profile in the two phases at the final time step. This time, only the left part of the domain is shown for a better visualization. The resulting temperature profile is almost constant in the liquid region, while we can observe the ice heating up in the solid region, as temperature increases from the ambient value when we get closer to the PCI.

\begin{table}
\begin{center}
\begin{tabular}{lcc} 
   \toprule
   \textbf{Parameter} & \textbf{phase 1} & \textbf{phase 2} \\ 
   \midrule
   $\rho$ $[\si{\kg\per\m\cubed}]$ & 999.88 & 916.8 \\
   $c_p$ $[\si{\joule\per\kg\per\kelvin}]$ & 4179.6 & 2090.0 \\
   $\kappa$ $[\si{\watt\per\m\per\kelvin}]$ & 0.5557 & 2.220 \\
   $\mu$ $[\si{\kg\per\m\per\s}]$ & 1.787 & 1e4 \\
   \midrule
   $h_m$ $[\si{\joule\per\kg}]$ & 333700 & - \\
   $T_m$ $[\si{\kelvin}]$ & 273 & - \\
   \bottomrule
\end{tabular}
\end{center}
\caption{Phase-change coupled 2D corner flow. The material properties for the two phases are shown.}
\label{table:2dProbeMaterials}
\end{table}

\subsection{Phase-change coupled 2D corner flow - effect of the outflow channel thickness}

In the previous sections we have shown the potential of our method to handle convection-coupled phase-change problems. We have also shown that the embedded phase-change interface resembles a no-slip boundary for large values of viscosity in the solid, see figure \ref{fig:lidDrivenComparisonU}. Our method hence covers physical regimes that imply a considerable evolution of the PCI as well as regimes that result in a more or less stationary phase boundary. As our computational approach allows to vary the scenario's geometry, we focus on multiple channel thicknesses for the last numerical example.

\begin{figure}

    \centering
  \begin{subfigure}[c]{0.49\textwidth}
    \includegraphics[trim={0cm 0cm 0cm 0cm},clip, width=\textwidth]{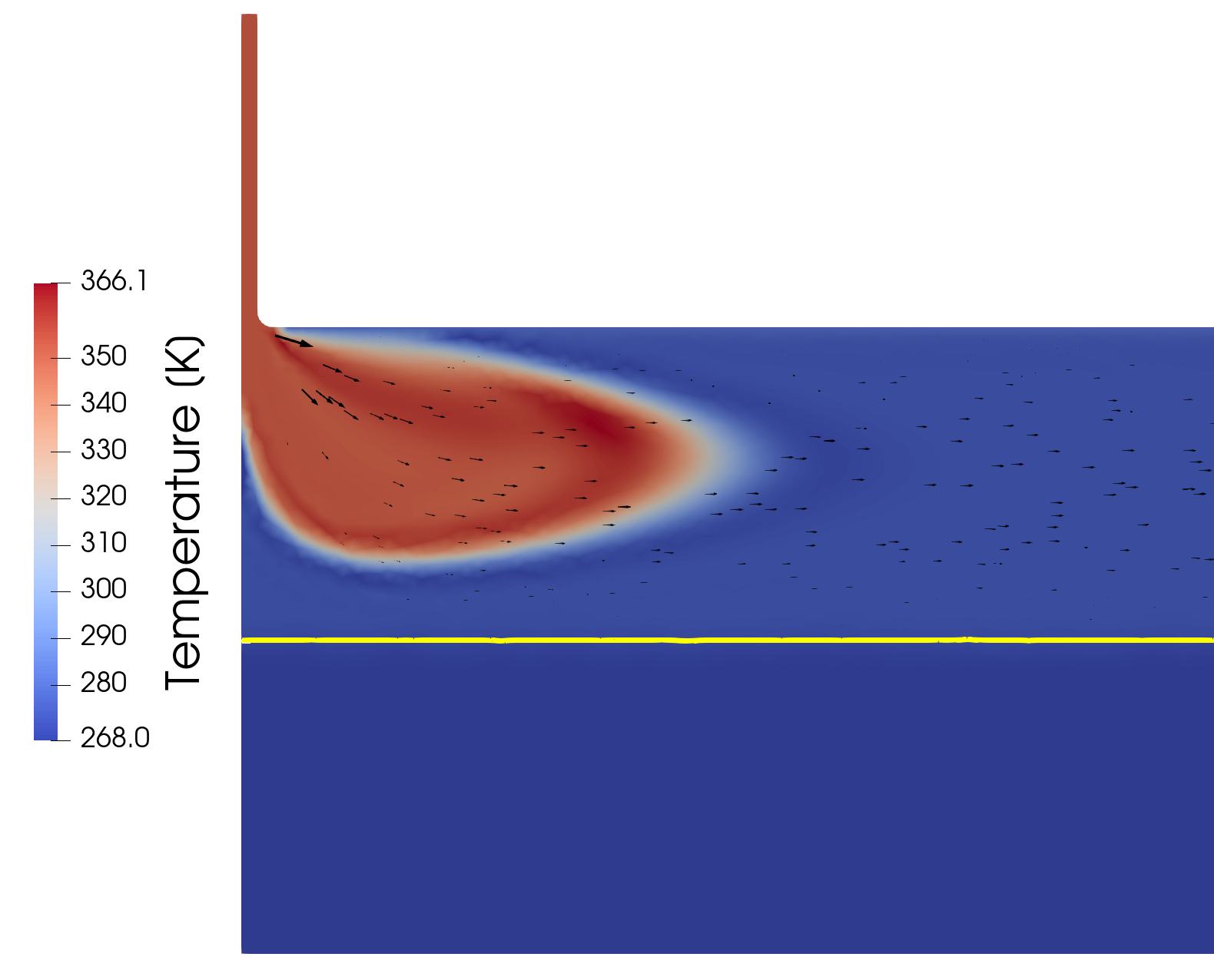}
    \caption{$t=10s$}
    \end{subfigure}
    \begin{subfigure}[c]{0.49\textwidth}
    \includegraphics[trim={0cm 0cm 0cm 0cm},clip, width=\textwidth]{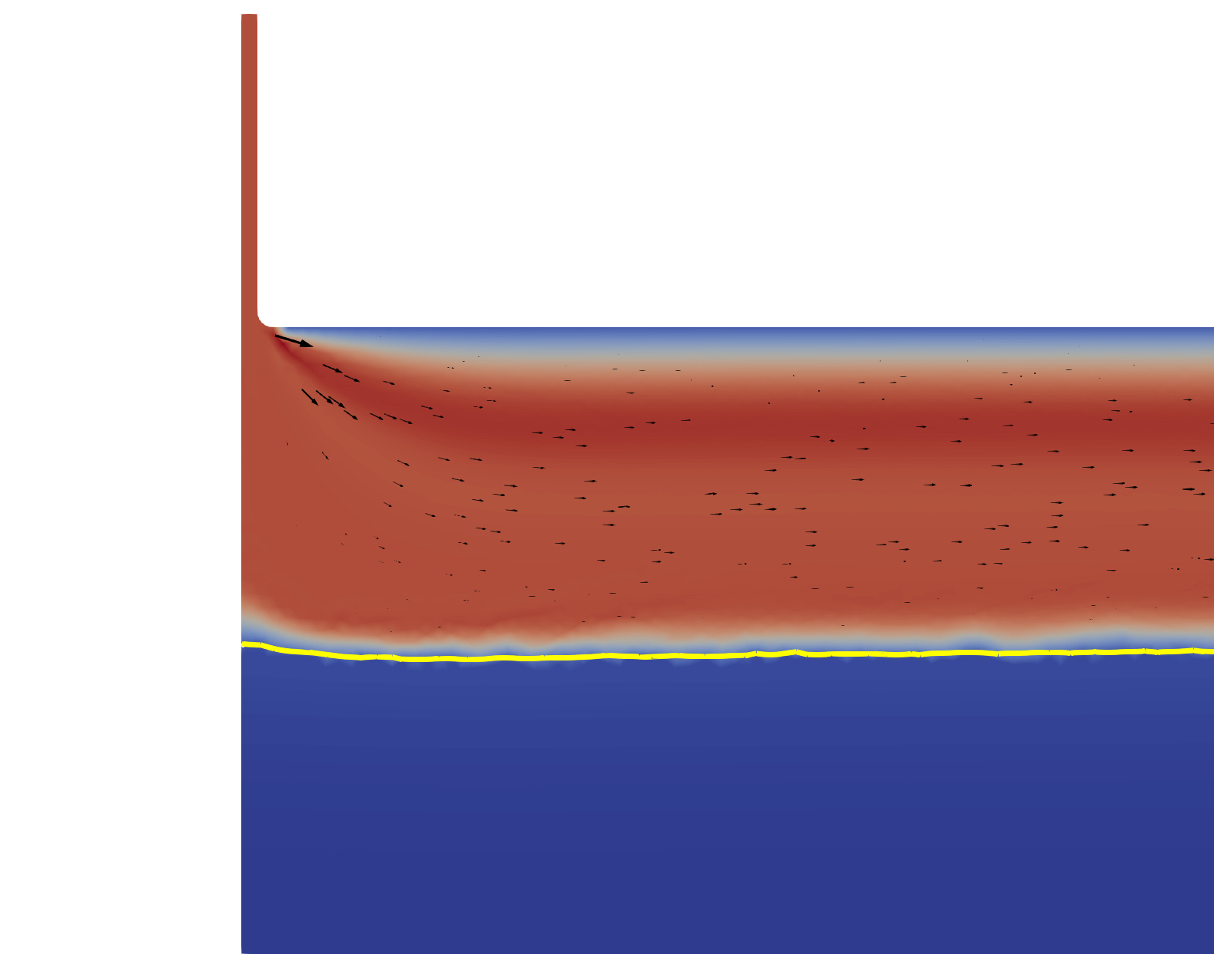}
    \caption{$t=2500s$}
    \end{subfigure}
    
    \begin{subfigure}[c]{0.49\textwidth}
    \includegraphics[trim={0cm 0cm 0cm 0cm},clip, width=\textwidth]{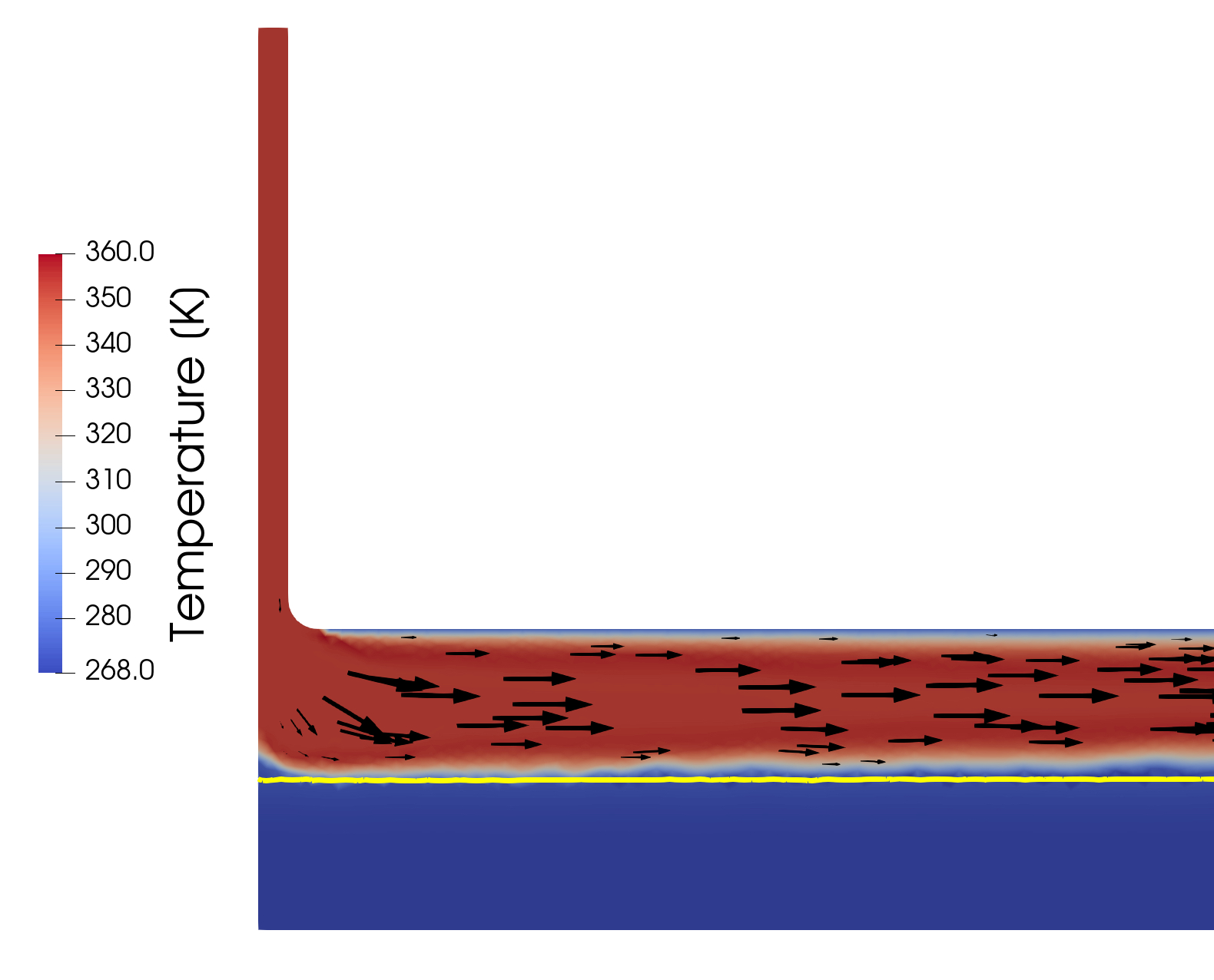}
    \caption{$t=10s$}
    \end{subfigure}
    \begin{subfigure}[c]{0.49\textwidth}
    \includegraphics[trim={0cm 0cm 0cm 0cm},clip, width=\textwidth]{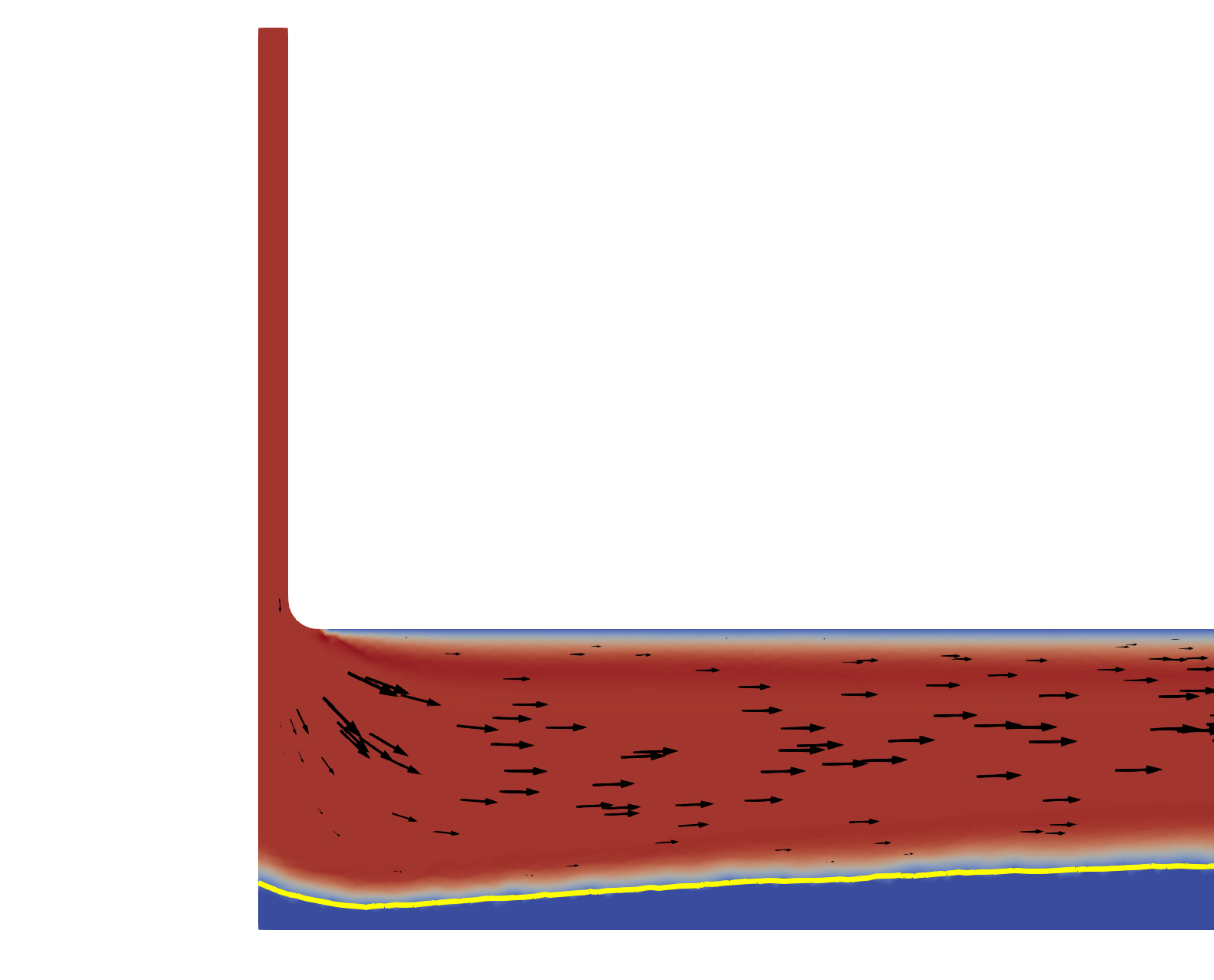}
    \caption{$t=2500s$}
    \end{subfigure}
    
    \caption{Phase-change coupled 2D corner flow - effect of the outflow channel thickness. We run the simulation of Fig.\ \ref{fig:2dProbeDomain} for two additional values of the channel thickness $d$. The temperature profile is plotted at $t=10s$ (left) and $t=2500s$ (right). The black arrows represent the velocity vectors at each point and their size is proportional to the second component of velocity. At the top, we show the thicker channel with $d=0.40$. At the bottom, we show the thinner channel with $d=0.10$. It is clearly visible that the thinner channel experiences a significant evolution of the PCI, whereas the PCI is almost stationary in the thicker channel.} 
    \label{fig:probe2DthicknessComparison}
    
\end{figure}

\begin{figure}

    \centering
    
    \includegraphics[trim={0cm 0cm 0cm 0cm},clip, width=0.45\textwidth]{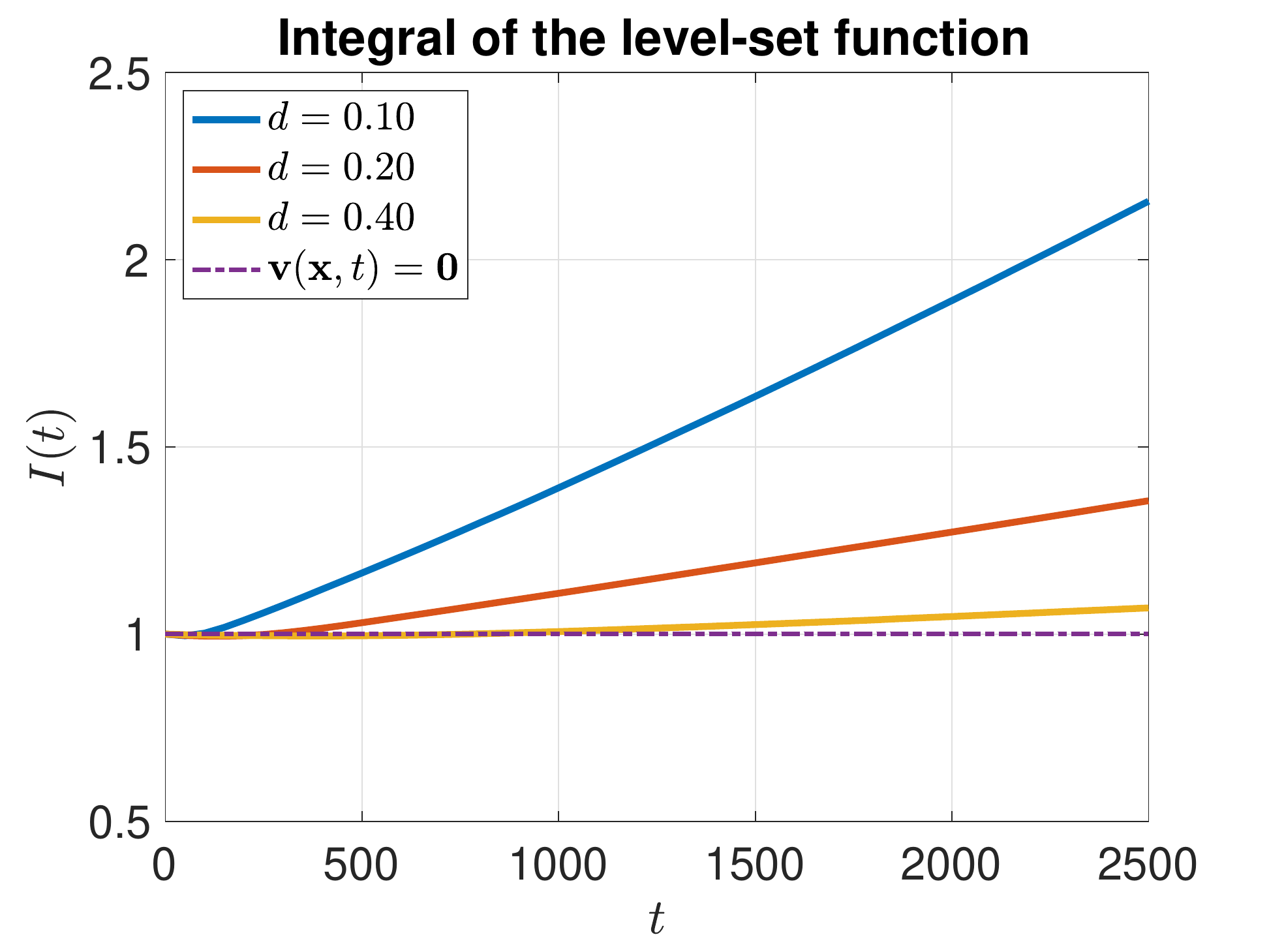}
    
    \caption{Phase-change coupled 2D corner flow - effect of the outflow channel thickness. The evolution of the liquid portion of the domain, computed as the integral of the level-set function on the liquid phase, is plotted over time for three different channel thickness values. The red line with $d=0.20$ refers to the simulation in Section \ref{2DmeltingProbe}. The purple dashed line denotes the reference case, in which the propagation velocity of the PCI is artificially set to zero to resemble the absence of phase change. The thinnest channel experiences the strongest increase in liquid area (blue), whereas the liquid area only moderately increases in the thickest channel (yellow).}
    \label{fig:probe2DthicknessMatlab}
    
\end{figure}

We consider the same setup as in Section \ref{2DmeltingProbe} and call $d$ the thickness of the outflow channel in Figure \ref{fig:2dProbeDomain}. In addition to the previous test case, we perform two simulations where we modify the thickness to $d=0.10$ and $d=0.40$, resulting in aspect ratios of $1/10$ and $1/40$. The initial position of the PCI is always located in the middle of the outflow channel, hence it shifts to $x=0.25$ and $x=0.40$, respectively. We introduce the same amount of heat into the system as before. Due to the varied outflow channel thickness and the updated position of the PCI, however, we expect to observe a different behaviour of the interface evolution. Figure \ref{fig:probe2DthicknessComparison} shows the temperature profile on both domains at $t=10s$ and $t=2500s$. As anticipated, the melting of the ice has a much larger effect on the thinner channel, to the point that the liquid covers the majority of the domain at the end of the simulation. In order to show the effect of the outflow channel thickness on the melting efficiency in a quantitative manner, we determine the amount of liquid in the channel by computing the integral of the level-set function over time in the liquid region, that is 

\begin{equation}
	I(t) = \frac{\int_{\Omega_1(t)}\Phi(t)\diff\Omega}{\int_{\Omega_{1,0}}\Phi_0\diff\Omega}.
\end{equation}
Here, the subscript $0$ denotes the initial time. Figure \ref{fig:probe2DthicknessMatlab} shows the plot of the liquid area $I(t)$ over time for the three different channel thicknesses that we have considered. For comparison, the purple dashed line represents a simulation without phase-change, i.e.\ the velocity of the PCI is artificially set to zero. We can clearly see that, depending on the properties of the domain, melting can deeply affect the outcome of the simulation. 

%% file: section5.tex
\section{Conclusion}
\label{section5}

In this work we presented a computational strategy for the numerical solution of phase-change problems with space-time finite elements. Two coupled problems for the flow field and temperature are considered, while the phase-change interface is tracked with a level-set method. The propagation velocity of the interface is determined by the Stefan condition, a heat-flux jump condition that accounts for local energy conservation across the phase-change interface. The Stefan condition constitutes the need for a numerical approximation of the heat flux at both sides of the interface. We leveraged the ghost-cell method, which was extended to the finite element framework. This technique considers two separate temperature problems for each phase, such that the melting temperature is enforced at the nodes close to the interface. Thus, the discontinuity of the heat flux at the interface can be retrieved. The ghost split does not require an enrichment of the finite element basis functions, like other methods. This simplifies the implementation in numerical codes.

We verified the 2D algorithm against a quasi 1D, single-phase Stefan problem. We showed convergence to the analytical solution both for structured quadrilateral and unstructured triangular meshes. Then we investigated the lid-driven cavity problem with phase change at its bottom boundary. Despite not imposing the no-slip condition directly at the interface, we obtained the same velocity profile as in a flow field simulation in the sole liquid region. Thus the implicit handling of the phase-change boundary reduces to a no-slip boundary as the viscosity in the solid phase tends to large values. Our method is hence applicable to multi-regime situations. Finally, we demonstrated the feasibility of our method in a more complex geometry. Inspired by one of our applications, we considered a 2D corner flow setup, in which warm water flowing through a thin inflow channel is diverted into a thicker outflow channel with a PCI at its outer side. This way we demonstrated the capability of our method to handle complex problems on domains of interest. The last simulation showed a quantitative analysis of the melting efficiency, where the liquid volume fraction increases more rapidly in thinner channels.

This framework will be the starting point to study more complex applications in 2D and 3D, and to guide computational model setups in the future. 